\documentclass[oneside]{amsart}% Math and Physical Sciences Numbered Reference Style

\usepackage{graphicx}%
\usepackage{multirow}%
\usepackage{amsmath,amssymb,amsfonts,bm}%
\usepackage{amsthm}%
\usepackage{mathrsfs}%
\usepackage[title]{appendix}%
\usepackage{xcolor}%
\usepackage{textcomp}%
\usepackage{manyfoot}%
\usepackage{booktabs}%
\usepackage{algorithm}%
\usepackage{algorithmicx}%
\usepackage{algpseudocode}%
\usepackage{listings}%
\usepackage{color}
\usepackage{geometry}
\usepackage{tabto}
\usepackage{calligra}
\usepackage{bbm}
\usepackage{cancel}
\usepackage{mathtools}
\usepackage{version}
\usepackage[math]{cellspace}
\usepackage{subfig}
\usepackage{tabularx}
\usepackage{arydshln}
\usepackage[abs]{overpic}
\usepackage{lineno}
\theoremstyle{thmstyleone}%
\theoremstyle{thmstyletwo}%

\theoremstyle{thmstylethree}%

\newcommand\red{\color{black}}
\newcommand\black{\color{black}}

\long\def\NOTE#1{}%{{\par\small\color{red}  {\noindent{\bf Nota}.~ #1 }\par\medskip}} %%

\newcommand{\Gp}{G_\mathrm{p}}
\newcommand{\Gt}{G_\mathrm{t}}
\newcommand{\Uii}{U_\mathrm{ii}}
\newcommand{\Uid}{U_\mathrm{id}}
\newcommand{\EGP}{EGP}
\newcommand{\RA}{Ra}
\newcommand{\E}{E}
\newcommand{\VG}{V_G}
\newcommand{\G}{G}
\newcommand\kuno{k_1}
\newcommand\kdue{k_2}
\newcommand{\U}{U}
\newcommand{\Il}{I_\ell}
\newcommand\Ip{I_p}
\newcommand\I{I}
\newcommand\VI{V_I}
\newcommand\muno{m_1}
\newcommand\mdue{m_2}
\newcommand\mtre{m_3}
\newcommand\mquattro{m_4}
\newcommand\mcinque{m_5}
\newcommand\msei{m_6}
\newcommand\HE{HE}
\newcommand\Se{S}
\newcommand\keuno{k_{e_1}}
\newcommand\kedue{k_{e_2}}
\newcommand\Fcns{F_{\mathrm{cns}}}
\newcommand\Vm{V_m}
\newcommand\Km{K_m}
\newcommand\X{X}
\newcommand\pdueU{p_{2U}}
\newcommand\Qsto{Q_\mathrm{sto}}
\newcommand\Qstouno{\Qsto^{(1)}}
\newcommand\Qstodue{\Qsto^{(2)}}
\newcommand\Qgut{Q_\mathrm{gut}}
\newcommand{\dotQ}[1]{\dot{Q}_\mathrm{#1}}
\newcommand\Kempt{k_\mathrm{empt}}
\newcommand\Kgri{k_\mathrm{gri}}
\newcommand\Kabs{k_\mathrm{abs}}
\newcommand\BW{BW}
\newcommand\D{D}
\newcommand\impulse{\delta}
\renewcommand\b{b}
\renewcommand\d{d}
\renewcommand\a{a}
\renewcommand\c{c}
\newcommand{\Ipo}{I_\mathrm{po}}
\newcommand\K{K}
\newcommand\h{h}
\newcommand{\Spo}{S_\mathrm{po}}
\newcommand{\Idot}[1]{\dot{I}_{\mathrm{#1}}}
\newcommand\Y{Y}
\newcommand\Iuno{I_1}
\newcommand\Id{I_\mathrm{d}}
\newcommand\ki{k_i}
\newcommand\kpuno{k_\mathrm{p1}}
\newcommand\kpdue{k_\mathrm{p2}}
\newcommand\kptre{k_\mathrm{p3}}
\newcommand\kpquattro{k_\mathrm{p4}}

\newcommand\Kmax{k_\mathrm{max}}
\newcommand\Kmin{k_\mathrm{min}}

\newcommand\Vmx{V_{mx}}
\newcommand\Vmzero{V_{m0}}
\newcommand\Kmzero{K_{m0}}

\newcommand\ttheta{\bm{\vartheta}}
\newcommand\loss{\mathcal{L}}
\newcommand\dt{\Delta t}
\newcommand\tol{\mathsf{tol}}
\newcommand\ind[1]{\mathbb{I}_{\{#1\}}}

\newcommand{\Gmax}{\overline{\G}}
\newcommand{\tmax}{\overline{t}}
\newcommand{\Gbio}{\G^{bio}}
\newcommand{\tbio}{t^{bio}}
\newcommand{\tolG}{\tol_{\G}}

\newcommand{ \rv }[1]{\textcolor{black}{#1}}
\newcommand\fc{\color{black}}
\newcommand\cf{\color{black}}
\usepackage{todonotes}
\date{}

\begin{document}
	
	\title[Automatic computation of the glycemic index: the glucose standard]{Automatic computation of the glycemic index:\\data driven analysis of the glucose standard}

	\author{Fabio Credali $^\dagger$}
	\address{Computer, Electrical and Mathematical Sciences and Engineering division, King Abdullah University of Science and Technology, Thuwal 23955, Saudi Arabia}
	
	\email{fabio.credali@kaust.edu.sa}
	\author{Maria Teresa Venuti $^\dagger$}
	\address{Dipartimento di Biologia e Biotecnologie ``Lazzaro Spallanzani'', Universit\`a degli Studi di Pavia, via Ferrata 9, 27100, Pavia, Italy}
	\thanks{$^\dagger$ These authors contributed equally to this work.}
	\email{mariateresa.venuti01@universitadipavia.it}
	\author{Daniele Boffi}
	\address{Computer, Electrical and Mathematical Sciences and Engineering division, King Abdullah University of Science and Technology, Thuwal 23955, Saudi Arabia; Dipartimento di Matematica ``F. Casorati'', Universit\`a degli Studi di Pavia, via Ferrata 5, 27100, Pavia, Italy}
	\email{daniele.boffi@kaust.edu.sa}
	\author{Paola Rossi}
	\address{Dipartimento di Biologia e Biotecnologie ``Lazzaro Spallanzani'', Universit\`a degli Studi di Pavia, via Ferrata 9, 27100, Pavia, Italy}
	\email{paola.rossi@unipv.it}
	%

%\begin{frontmatter}
	
%	\title{Automatic computation of the glycemic index:\\data driven analysis of the glucose standard}
%
%    \author[fabio]{Fabio Credali\fnref{equal,corr}}
%	\ead{fabio.credali@kaust.edu.sa}
%	
%	\author[mteresa]{Maria Teresa Venuti\fnref{equal}}
%	\ead{mariateresa.venuti01@universitadipavia.it}
%	
%	\author[fabio,daniele]{Daniele Boffi}
%	\ead{daniele.boffi@kaust.edu.sa}
%	
%	\author[mteresa]{Paola Rossi}
%	\ead{paola.rossi@unipv.it}
%	
%	\address[fabio]{CEMSE Division, King Abdullah University of Science and Technology, Thuwal, 23955, Saudi Arabia}
%	\address[mteresa]{Dipartimento di Biologia e Biotecnologie ``Lazzaro Spallanzani'', Universit\`a degli Studi di Pavia, via Ferrata 5, Pavia, 27100, Italy}
%    \address[daniele]{Dipartimento di Matematica ``F. Casorati'', Universit\`a degli Studi di Pavia, via Ferrata 9, Pavia, 27100, Italy}
%	\fntext[equal]{These authors contributed equally to this work.}
%	\fntext[corr]{Corresponding author}
	
	\begin{abstract}
		The Glycemic Index (GI) is a tool for classifying carbohydrates based on their impact on postprandial glycemia, useful for diabetes prevention and management. This study applies a mathematical model for a data driven simulation of the glycemic response following glucose ingestion. \rv{The analysis is performed on a dataset of 35 healthy subjects undergone a standard 50 g oral glucose test. The results reveal} a direct correlation between glucose response profiles and parameters describing glucose absorption, enabling the classification of subjects into three groups based on the timing of their glycemic peak: \rv{$<30$ min, $30-50$ min, $>50$ min}. \rv{These findings highlight the ability of a physiology-based mathematical model to capture inter-individual variability in postprandial glucose dynamics and represent a step toward simulation-based approaches for GI estimation.}\medskip
        %Our results offer potential applications for both glycemic index simulation and advancing biological studies on diabetes.\medskip
		
		\medskip
		
	\end{abstract}
	
%\end{frontmatter}

\maketitle

\section{Introduction}

Diabetes is a major global health concern, with 589 million adults (20-79 years) affected and 3.4 million deaths reported in 2024 by the International Diabetes Federation (IDF). The total number of people living with diabetes is predicted to rise to 853 million by 2050 ~\cite{IDF2025}.
The $90\%$ of these cases are type~2 diabetes (T2DM), making it the predominant form of the disease worldwide ~\cite{IHME2024,2IHME2024}.

T2DM is characterized by peripheral insulin resistance and progressive pancreatic $\beta-$cell dysfunction, resulting in chronic hyperglycemia. Due to its insidious and often asymptomatic onset, T2DM is frequently diagnosed at a late stage, increasing the likelihood of complications at the time of diagnosis ~\cite{Gregg2014,King1999}. The etiology of T2DM is multifactorial, involving a combination of genetic predisposition and environmental factors. Major risk factors include excess body weight, advancing age, certain ethnic backgrounds, and a family history of diabetes.
Effective management relies on comprehensive lifestyle modification, emphasizing a balanced diet, regular physical activity, weight management, and smoking cessation. When lifestyle changes alone are insufficient to maintain glycemic control, pharmacological intervention is introduced, with metformin serving as the standard first-line therapy ~\cite{IDF2025}.

Lifestyle interventions are effective in preventing or delaying T2DM progression~\cite{uusitupa2019prevention,kondrad2023metformin,ramachandran2006indian}. Scientific consensus highlights the impact of dietary choices on T2DM risk, with high glycemic index/load (GI/GL) diets potentially contributing to disease onset, though causality remains uncertain. The concept of GI was introduced in 1981 to classify carbohydrates based on their effect on postprandial glycemia \cite{jenkins1981glycemic}. It represents the blood glucose response of a 50-gram carbohydrate portion of food, expressed as a percentage of the same amount of carbohydrate from a reference food (usually pure glucose). The GI ranks the glycemic potential of carbohydrates in different food. Various factors influence GI, including carbohydrate type, starch properties, food processing, and macronutrient interactions~\cite{jenkins1981glycemic,russell2016impact,lovegrove2017role}.                                                                      
Given the increasing attention on the glycemic index (GI) of food, the International Organization for Standardization (ISO) has established the official standard for measuring the GI~\cite{iso}. This method involves giving ten or more healthy individuals a portion of food containing 50 grams of digestible carbohydrates, measuring blood glucose levels before eating and at regular intervals for two hours after eating, and plotting the changes in blood glucose concentration over time as a curve. %The GI is then calculated as the incremental area under the glucose curve after eating the test food, divided by the corresponding area after consuming a control food, and multiplying the result by 100 to obtain a percentage.
The GI is then calculated as the incremental area under the glucose curve after eating the test food expressed as a percentage of the corresponding area after consuming the control food.
The ISO method defines the GI, outlines qualifying factors, specifies requirements for its application, and recommends criteria for classifying food into low, medium, and high GI. It has been calibrated by independent laboratories. The GI serves as a tool for comparing and understanding the biological effects of different carbohydrates.
In addition, participants should fast for at least 12 hours, avoid meals rich in carbohydrates and fats, refrain from alcohol and smoking, avoid intense physical exercise in the days prior, and ensure adequate sleep to minimize confounding factors and ensure accurate test results. ~\cite{iso,Robertson2002,siler1998inhibition,malkova2000prior}

Accurately determining GI is crucial for dietary planning and consumer guidance. Traditional GI testing involves controlled human trials, but individual variability and small sample sizes limit precision. Advancements in modeling offer promising alternatives for real-time glucose monitoring and improved GI estimation, benefiting both consumers and the food industry. More precise GI measurement methods would enable manufacturers to provide reliable GI values on food labels, supporting better dietary choices, diabetes management, and public health initiatives~\cite{brand2020relationship,marsh2011glycemic,livesey2019dietary,wolever2013glycaemic,barclay2021dietary}.

During the last decades several mathematical models for the glucose-insulin system have been proposed to provide effective simulation methods for treating and preventing diabetes and obesity. Such models are usually classified into two families. \textit{Minimal models}~\cite{bergman79,breda01,dallaman02,campioni09} just describe the key components of a system, allowing the measurement of non-accessible data such as insulin sensitivity and $\beta-$cells responsivity. On the other hand, \textit{maximal models} give a comprehensive description on the ingestion process~\cite{hovorka04,dallaman07,degaetano13,moxon16,viceconti17,visentin19}. The main features of several minimal and maximal models are summarized in~\cite{cobelli14,dallaman22}. \rv{All these models are based on ordinary differential equations (ODEs) only describing the temporal variation of the quantities of interest: the spatial homogeneity of compartments is an admissible assumptions when analyzing the glucose-insulin system from a macroscopic point of view. On the other hand, the inclusion of diffusion models based on partial differential equations (PDEs) is necessary when studying local phenomena, see e.g.~\cite{buchwald,aguilar}.}

In this work we consider the Dalla Man--Rizza--Cobelli maximal model. This model was introduced in~\cite{dallaman07} as a first example of \textit{meal simulator} and consists of six compartments \rv{of ODEs} representing glucose and insulin kinetics, gastro-intestinal tract, liver, kidneys and adipose tissues. It describes the digestion process from different points of view beyond just intravenous glucose perturbations~\cite{cobelli82,cobelli83,vicini99}. The Dalla Man--Rizza--Cobelli model represents a milestone in the field of bioengineering for diabetes prevention since it has been applied several times and for different purposes (see e.g.~\cite{visentin18,sanchez19,jacobs23}). As a consequence, several features of the model have been addressed during the past years~\cite{dallaman06,toffolo06,dalla07gim}, as well as enhancements and improvements~\cite{piccinini15}.

\fc
Despite the availability of standardized ISO procedures, GI determination still relies on repeated \textit{in vivo} testing on healthy subjects, which is time-consuming, costly, and inherently affected by inter-individual variability. A quantitative and reproducible simulation framework capable of capturing postprandial glucose dynamics could represent a valuable complementary tool to experimental protocols.
In this work, we propose a data-driven modeling approach based on the Dalla Man–Rizza–Cobelli model to simulate the glycemic response to oral glucose ingestion. The study is conducted exclusively on healthy subjects, consistently with the ISO definition of GI.
The objectives of this study are:
\begin{enumerate}%[i)]
    \item to calibrate the model using \textit{in vivo} glucose measurements obtained after a standard 50 g oral glucose test;
    \item to investigate the relationship between glucose absorption parameters and the timing of the glycemic peak;
    \item to evaluate whether model-based parameter identification enables classification of distinct metabolic response profiles.
\end{enumerate}
This work represents a first step toward the development of a simulation-based framework for GI estimation, potentially reducing reliance on extensive human testing.
\cf

%Since our aim is not the direct study of diabetes, but rather the design of an effective simulation model for the glycemic index, we consider the basic formulation introduced in~\cite{dallaman07}. The aim of the present work is to apply the model to simulate the glycemic response after ingestion of pure glucose, that is the reference for the GI computation. We apply the model to a dataset of \textit{in vivo} measurements on healthy subjects to fit the data and estimate the physiological parameters governing the system. Our analysis describes the direct correlation between the profile of the glucose response and the model parameters describing glucose absorption, providing a useful tool for both GI simulation and biological studies on diabetes.

This paper is organized as follows. We recall the main features of the Dalla Man--Rizza--Cobelli model in Section~\ref{sec:system}. In Section~\ref{sec:optimization} we describe the optimization process we applied for the parameter estimation. The results and the related discussion are then collected in Section~\ref{sec:results}. Finally, in Section~\ref{sec:conclusions}, we draw some conclusions and mention future advancements.

\section{Model of the glucose--insulin system}\label{sec:system}

We consider the Dalla Man--Rizza--Cobelli model of the glucose-insulin system introduced in~\cite{dallaman07}. The biological system is described by several ordinary differential equations grouped into six compartments, representing the process of digestion:
\fc
\begin{itemize}
    \item the \textit{gastro-intestinal tract} describes transit and absorption of glucose in stomach and intestine,
    \item the \textit{glucose system} describes the dynamics of glucose in plasma and tissues,
    \item \textit{muscle and adipose tissues} consume the absorbed glucose,
    \item the \textit{liver} equilibrates the endogenous glucose production depending on the plasma glucose and insulin,
    \item $\bm{\beta}$\textit{-cells} are responsible for the production of new insulin after glucose ingestion,
    \item the \textit{insulin system} describes the dynamics of insulin in plasma and liver.
\end{itemize}
\cf

\begin{figure}[t]
	\centering
	\includegraphics[width=0.9\linewidth]{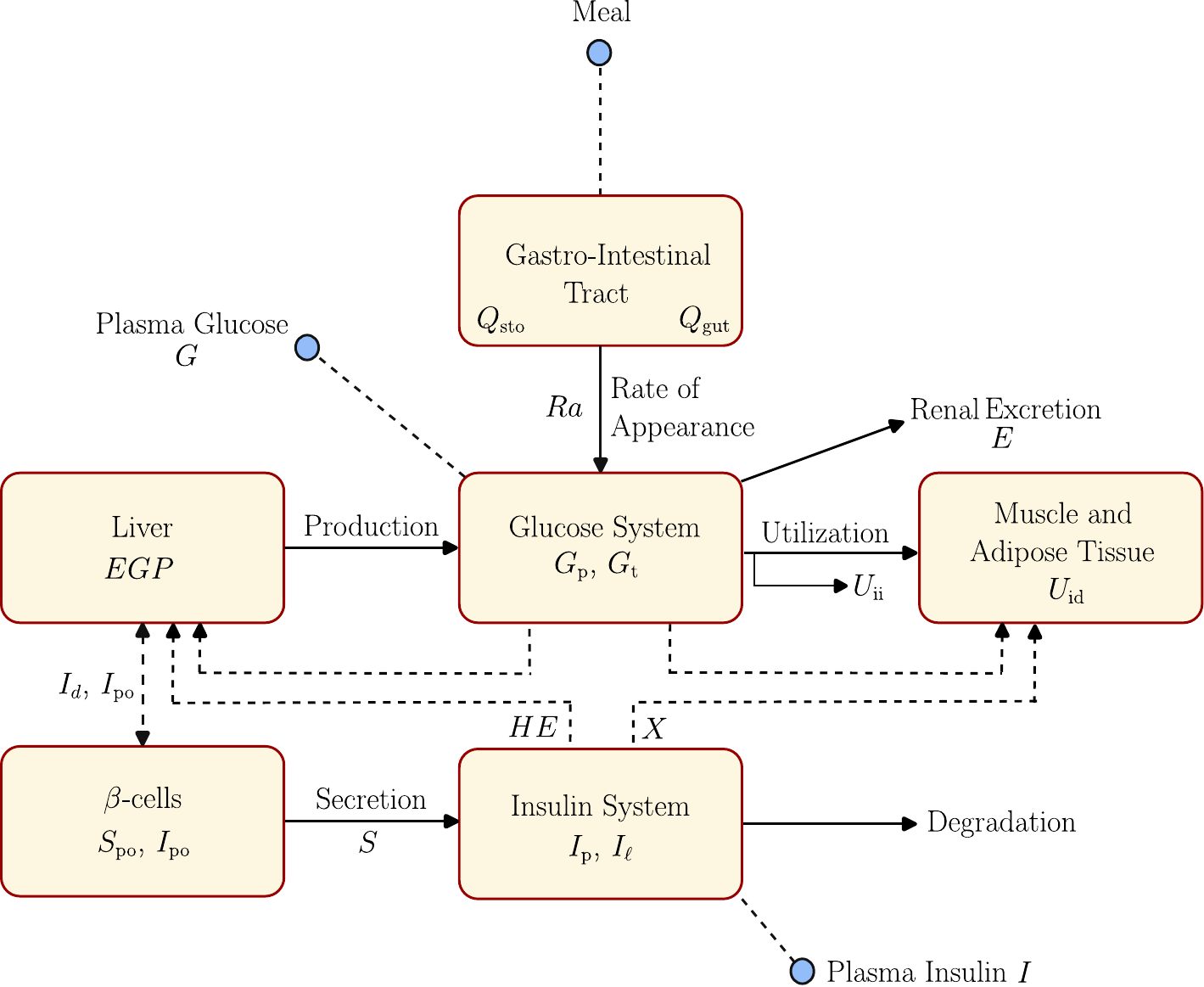}
	\caption{Representation of the Dalla Man--Rizza--Cobelli model with six compartments (yellow boxes) and their relations.}
	\label{fig:model}
\end{figure}

%Such compartments are the gastrointestinal tract, the liver, the glucose absorption system, the muscle and adipose tissues, the pancreas ($\beta-$cells), and the insulin kinetics.
A schematic representation of the compartments, together with their interplay, is shown in Figure~\ref{fig:model}. \rv{Moreover, Table~\ref{tab:compartments} collects all the variables and their meaning, while Tables~\ref{tab:intervals}-\ref{tab:param_dallaman} collect the model parameters.}
We now recall all the governing equations.% and the meaning of all involved variables and biological parameters \rv{(see also Tables~\ref{tab:compartments}-\ref{tab:param_dallaman})}.

\begin{table}\renewcommand{\arraystretch}{1.2}
    \centering
    \color{black}
    \begin{tabular}{lll}
    
        \multicolumn{3}{c}{\textbf{Summary of compartments}}\\
         %\textbf{Symbol}&\textbf{Unit}&\textbf{Meaning}\\
         \hline
         \multicolumn{3}{l}{\textit{Gastrointestinal tract}}\\
         \hline
         $\Qstouno$ & mg & glucose mass in stomach (solid phase)\\
         $\Qstodue$ & mg & glucose mass in stomach (liquid phase)\\
         $\Qsto$ & mg & total glucose mass in stomach\\
         $\Qgut$ & mg & glucose mass in intestine\\
         \hline
         \multicolumn{3}{l}{\textit{Glucose dynamics (Glucose subsystem + Liver + Utilization)}}\\
         \hline
         $\Gp$ & mg/kg & glucose mass in plasma rapidly equilibrating tissues\\
         $\Gt$ & mg/kg & glucose mass in slowly equilibrating tissues\\
         $\G$ & mg/dL & plasma glucose concentration\\
         $\RA$ & mg/kg/min & rate of appearance of glucose in plasma\\
         $\EGP$ & mg/kg/min & endogenous glucose production\\
         $\Uii$ & mg/kg/min & glucose utilization (insulin independent)\\
         $\Uid$ & mg/kg/min & glucose utilization (insulin dependent)\\
         $\U$ & mg/kg/min & total glucose utilization\\
         $\E$ & mg/kg/min & renal excretion\\
         \hline
         \multicolumn{3}{l}{\textit{Insulin dynamics (Insulin subsystem + $\beta$-cells)}}\\
         \hline
         $\Ip$ & pmol/kg  & insulin mass in plasma\\
         $\Il$ & pmol/kg & insulin mass in liver\\
         $\I$ & pmol/L & plasma insulin concentration\\
         $\Ipo$ & pmol/kg & insulin in portal vein\\
         $\Id$ & pmol/L & delayed insulin signal\\
         $\HE$ & dimensionless & hepatic insulin extraction\\
         $\X$ & pmol/L & insulin in interstitial fluid\\
         $\Se$ & pmol/kg/min & insulin secretion (portal vein to liver)\\
         $\Spo$ & pmol/kg/min & insulin secretion (pancreas to portal vein)\\
        \hline
    \end{tabular}
    \color{black}
    \caption{Main model compartments. A summary of the related parameters is reported in Table~\ref{tab:intervals} and~\ref{tab:param_dallaman}.}
    \label{tab:compartments}
\end{table}

\subsection{Glucose subsystem}

The glucose subsystem consists of two equations describing the evolution of glucose masses in plasma and rapidly equilibrating tissues $\Gp$ (mg/kg) and in slowly equilibrating tissues $\Gt$ (mg/kg)
\begin{equation}\label{eq:glucose}
	\left\{
	\begin{aligned}
		&\dot{\Gp}(t) = \EGP(t)+\RA(t)-\Uii(t)-\E(t)-\kuno\, \Gp(t)+\kdue\,\Gt(t)&&\qquad\Gp(0) = \G_{\mathrm{p}b}\\
		&\dot{\Gt}(t) = -\Uid(t)+\kuno\, \Gp(t) -\kdue \,\Gt(t)&&\qquad\Gt(0) = \G_{\mathrm{t}b}\\
		&\G(t) = \frac{\Gp}{\VG}.&&
	\end{aligned}
	\right.
\end{equation}

The plasma glucose concentration $\G$ (mg/dL) is obtained dividing $\Gp$ by the distribution volume of glucose $\VG$ (dL/kg), \rv{that is the factor describing how glucose distributes into plasma}. In the first equation, $\EGP$ (mg/kg/min) is the endogenous glucose production, $\RA$ (mg/kg/min) denotes the rate of appearance of glucose in plasma, and $\E$ (mg/kg/min) the renal excretion. The letter $\U$ (mg/kg/min) stands for the glucose utilization which is classified into insulin independent $\Uii$ and insulin dependent $\Uid$ (see Section~\ref{sec:utilization}). Finally, $\kuno$ and $\kdue$ (min$^{-1}$) are two rate parameters, while the initial conditions are the basal values $\G_{\mathrm{p}b}$, $\G_{\mathrm{t}b}$ of $\Gp$ and $\Gt$, respectively.

\subsection{Insulin subsystem}

The insulin subsystem is governed by two equations describing the relation between the insulin mass in plasma $\Ip$ (pmol/kg) and in the liver $\Il$ (pmol/kg). The equations read
\begin{equation}
	\left\{
	\begin{aligned}
		&\dot{\Il}(t) = -(\muno+\mtre(t))\, \Il(t) + \mdue\, \Ip(t)+\Se(t)&&\qquad\Il(0)=I_{\ell b}\\
		&\dot{\Ip}(t) = -(\muno+\mquattro)\, \Ip(t)+\muno\, \Il(t)&&\qquad\Ip(0)=I_{pb}\\
		&\I(t) = \frac{\Ip}{\VI}.&&\\
	\end{aligned}
	\right.
\end{equation}

The variable $\Se$ (pmol/kg/min) represents the insulin secretion, which is further described in Section~\ref{sec:insulin_secretion}. The insulin concentration in plasma $\I$ (pmol/L) is obtained by dividing $\Ip$ by the distribution volume of insulin $\VI$ (kg$^{-1}$). The initial conditions $I_{\ell b}$ and $I_{pb}$ are the basal values of $\Il$ and $\Ip$, respectively.

The coefficients $\muno,\,\mdue,\,\mquattro$ (min$^{-1}$) are fixed rate parameters. On the other hand, the coefficient $\mtre$ is expressed as a function of the hepatic insulin extraction $\HE$, we have
\begin{equation}
	\HE(t) = -\mcinque\,\Se(t)+\msei,
	\qquad
	\mtre(t) = \frac{\muno\,\HE(t)}{1-\HE(t)},
\end{equation}
where, $\mcinque$ and $\msei$ (min$^{-1}$) are constant rate parameters too. \rv{More precisely, $\msei$ represents the baseline extraction of insulin by the liver, which take place when insulin secretion is zero. Then, $\HE$ decreases as $\Se$ increases.}

\fc 
The coefficients $\muno,\dots,\msei$ play a crucial role in computing the basal values of $\Se$, $\Ip$ and $\Il$. Indeed, given the basal value $\HE_b$ of the hepatic insulin extraction, the following relations hold (see ~\cite[Sect. III-B]{dallaman07} for more details):
\begin{equation}
    \begin{aligned}
        \mtre(0) = \frac{\muno\,\HE_b}{1-\HE_b},
        \qquad
        \Se_b = \frac{\msei-\HE_b}{\mcinque},
        \qquad
        I_{pb} = \frac{\Se_b\,(1-\HE_b)}{\mdue\,\HE_b+\mquattro},
        \qquad
        I_{\ell b} = \frac{\Se_b-\mquattro\,I_{pb}}{\mtre(0)}.
    \end{aligned}    
\end{equation}
\cf 

%More information regarding $\HE$ and its relation with $\mtre$ at basal level can be found in~\cite[Sect. III-B]{dallaman07}. 

\subsection{Endogenous Glucose Production}

The Endogenous Glucose Production (EGP) is expressed in terms of glucose in plasma $\Gp$, portal insulin $\Ipo$ (pmol/kg) and delayed insulin signal $\Id$ (pmol/L), as
\begin{equation}\label{eq:egp}
	\begin{aligned}
		&\EGP(t) = \kpuno - \kpdue\, \Gp(t)-\kptre\, \Id(t)-\kpquattro\, \Ipo(t)\\
		&\EGP(0)=\EGP_b.
	\end{aligned}
\end{equation}
More precisely, $\Ipo$ is the amount of insulin in the portal vein and its associated equation will be described in Section~\ref{sec:insulin_secretion}. On the other hand, $\Id$ is governed by the following pair of equations

\begin{equation}\label{eq:id}
	\left\{
	\begin{aligned}
		&\Idot{1}(t) = -\ki\,(\Iuno(t)-\I(t))&&\qquad \Iuno(0) = \I_b\\
		&\Idot{d}(t) = -\ki\,(\Id(t)-\Iuno(t))&&\qquad \Id(0) = \I_b,
	\end{aligned}
	\right.
\end{equation}
where $\Iuno$ is an auxiliary variable. Regarding the involved parameters, $\kpuno$ (mg/kg/min) is the extrapolated $\EGP$ at zero glucose and insulin, indeed at basal state the following relation holds
\begin{equation}\label{eq:egpb}
	\kpuno = \EGP_b + \kpdue\,\G_{\mathrm{p}b} + \kptre\,\I_b + \kpquattro\,\I_{\mathrm{po},b}
\end{equation}
Furthermore, $\kpdue$ (min$^{-1}$) is the liver glucose effectiveness, $\kptre$ (mg/kg/min per pmol/l) and $\kpquattro$ (mg/kg/$\backslash$min/(pmol/kg))  govern the amplitude of the insulin action on the liver. Finally, $\ki$ (min$^{-1}$) takes into account the delay between insulin signal and insulin action. \red $\EGP$ denotes the glucose production carried out by the liver, thus only nonnegative values are admissible and a constraint is added to the model~\cite{rizza2016}. \black
%We may be led to think that $\EGP\ge0$ means production of glucose, while $\EGP<0$ represents its storage: this observation is wrong and $\EGP$ must be constrained to be nonnegative (the interested reader may refer to~\cite{rizza2016} for more details).

\subsection{Gastrointestinal tract}\label{sec:gastro_tract}

This set of equations describes the transit of glucose in stomach and intestine, as well as its absorption. The total amount of glucose is denoted by $\Qsto$ (mg), that is the sum of glucose in solid phase $\Qstouno$ and triturated phase $\Qstodue$. The glucose mass in the intestine is denoted by $\Qgut$ (mg). We have the following equations
\begin{equation*}
	\left\{
	\begin{aligned}
		&\dotQ{sto}^{(1)}(t) = -\Kgri\, \Qstouno(t) + \D\,\impulse(t)&&\quad \Qstouno(0)=0\\
		&\dotQ{sto}^{(2)}(t) = -\Kempt(\Qsto)\, \Qstodue(t)+\Kgri\, \Qstouno(t)&&\quad \Qstodue(0)=0\\
		&\dotQ{gut}(t) = -\Kabs\, \Qgut(t)+\Kempt(\Qsto)\, \Qstodue(t)&&\quad \Qgut(0)=0\\
		&\Qsto(t) = \Qstouno + \Qstodue&&\\
		&\RA(t) = \frac{f\, \Kabs\, \Qgut(t)}{\BW}.
	\end{aligned}
	\right.
\end{equation*}

The letter $\D$ (mg) denotes the quantity of ingested glucose, while $\impulse$ \red(dimensionless) \black is an impulse function (with Gaussian profile, for instance).
Again, $\RA$ is the glucose rate of appearance in plasma, depending on the body weight $\BW$ (kg) and on the fraction of intestinal absorption actually appearing in plasma. As initial conditions, the quantities of glucose in stomach and intestine are both assumed to be null. Finally, $\Kabs$ (min$^{-1}$) is the rate of intestinal absorption, $\Kgri$ (min$^{-1}$) the rate of grinding and $\Kempt(\Qsto)$ the rate of gastric emptying~\cite{dallaman06,dallaman07}.

\begin{figure}
	\centering
	\subfloat[]{\includegraphics[width=0.34\linewidth]{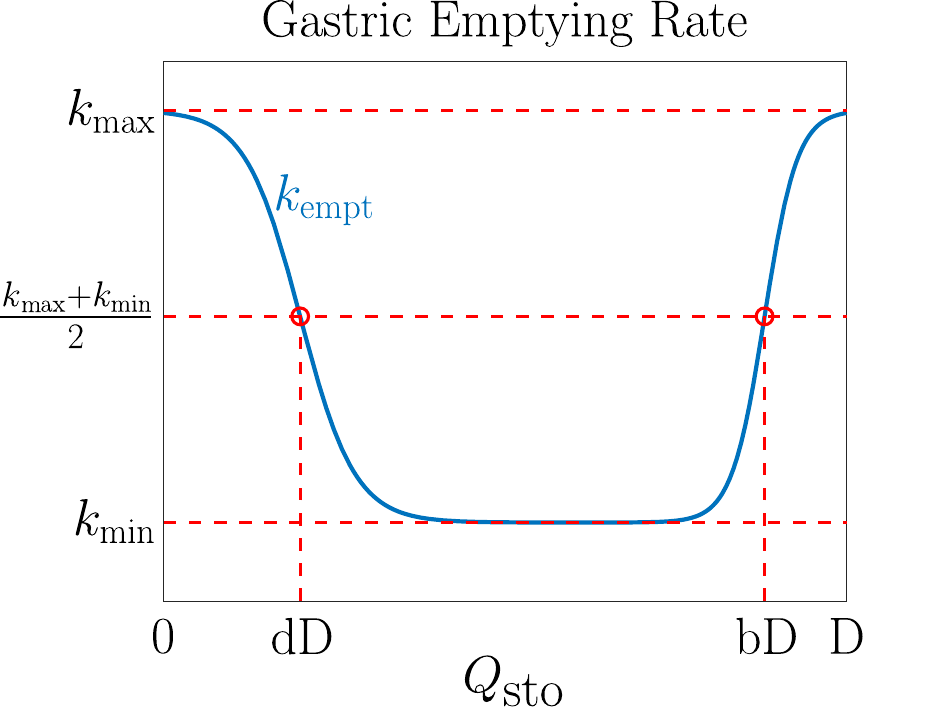}}
	\subfloat[]{\includegraphics[width=0.34\linewidth]{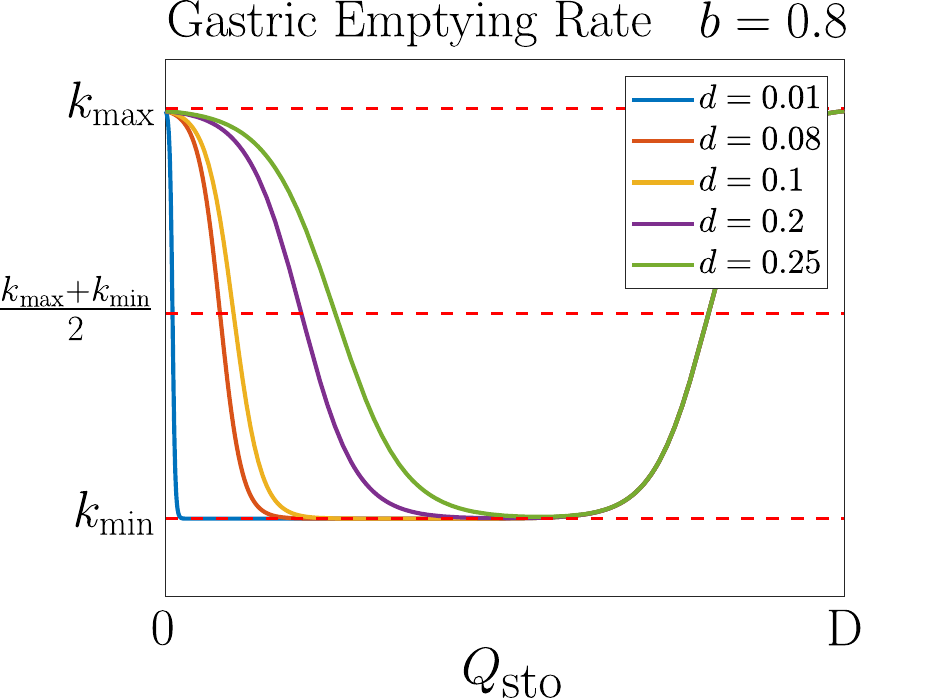}}
	\subfloat[]{\includegraphics[width=0.34\linewidth]{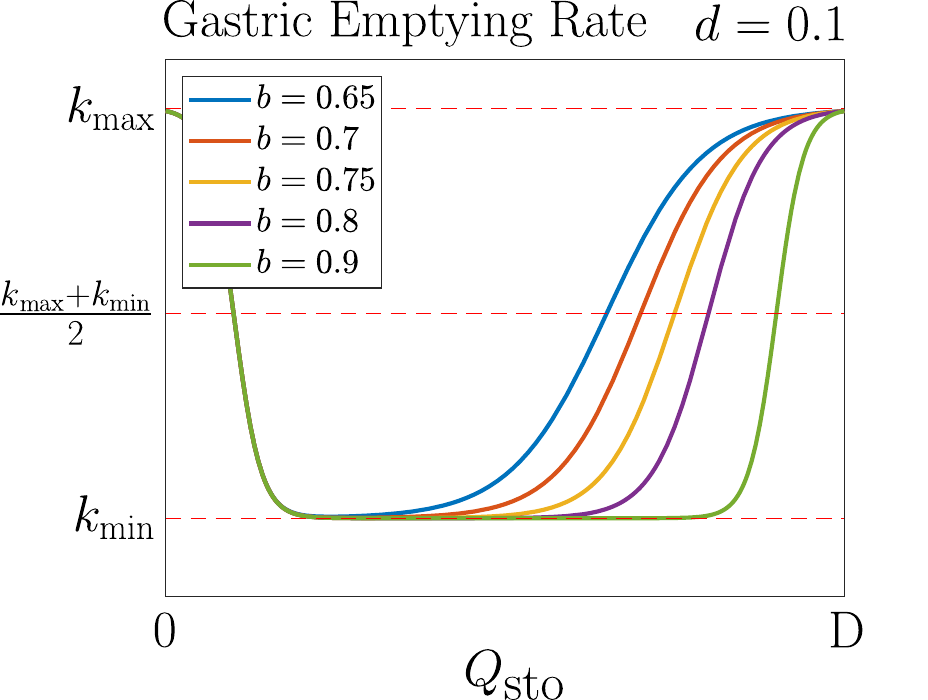}}
	\caption{\textsc{(a)} Profile of gastric emptying rate $\Kempt$~\rv{in Eq.~\eqref{eq:kempt}} for fixed values of $d$ and $b$ \rv{($0.2$ and $0.9$, respectively)}. \textsc{(b)} Profile of $\Kempt$ in variation of $d$ with $b=0.8$. \textsc{(c)} Profile of $\Kempt$ in variation of $b$ with $d=0.1$. \rv{The three plots were generated with fictitious values for illustration purposes. In particular, $D=10^4$, $\Kmax=0.05$, $\Kmin=0.008$.} }
	\label{fig:kempt}
\end{figure}

%We point out that $\Kempt$ depends on the quantity $\Qsto$ as prescribed by the following relation
\rv{The mathematical description of the gastric emptying rate was introduced in~\cite{dallaman06}. More precisely, }$\Kempt$ depends on the quantity $\Qsto$ as prescribed by the following relation

\begin{equation}\label{eq:kempt}
	\Kempt(\Qsto) = \Kmin+\left(\frac{\Kmax-\Kmin}2\right)
	(\tanh(\a(\Qsto-\b\D))
	-\tanh(\c(\Qsto-\d\D))+2).
\end{equation}
The profile of such function is shown in Figure~\ref{fig:kempt}. More precisely, $\Kempt=\Kmax$ when the stomach contains the amount $\D$ of ingested glucose. Then it decreases to its minimum $\Kmin$ during the digestion process and returns back to $\Kmax$ when the stomach is empty. The parameter $\b$ is the percentage of glucose dose for which $\Kempt$ decreases at the average $(\Kmax+\Kmin)/2$, as well as $\d$ is the percentage of dose for which $\Kempt$ is back to $(\Kmax+\Kmin)/2$. \rv{Notice that $\d\D$ corresponds to the inflection point in the decreasing phase, while $\b\D$ corresponds to the inflection point in the recovery phase.} 
\fc
In other words, the behavior of $\Kempt$ goes through the following phases:
\begin{enumerate}
    \item $\Kempt=\Kmax$ when the quantity of glucose in the stomach $\Qsto=0$, i.e. the stomach is ready for initiating the digestion process;
    \item for $\Qsto>0$ digestion starts and $\Kempt$ gradually reduces to $\Kmin$, $\d$ is the percentage of ingested glucose corresponding to the inflection point in this decay phase;
    \item $\Kempt=\Kmin$ during the core of the digestion process;
    \item when the digestion process is almost finished, $\Kempt$ smoothly recovers to $\Kmax$ and the curve has an inflection point in the correspondence of $\b\D$.
\end{enumerate}
\cf

The decay and recovery rates $\a$ and $\c$, respectively, depend on $\b$, $\d$ and $\D$. Indeed, at the time of ingestion, $\Kempt$ decreases with rate
$$
a = \frac{5}{2\D(1-\b)},
$$
whereas the recovery rate at the end of the digestion process is
$$
\c = \frac{5}{2\d\D}.
$$

We remark that if $\b$ is in a large percentage, then the decay rate $a$ is also large: this implies a slow digestion since $\Kempt$ rapidly reaches its minimum. A similar reasoning applies to $d$. If $\d$ is small, then the rate $\c$ is large, implying a fast recovery phase from $\Kmin$ back to $\Kmax$.

\subsection{Glucose utilization}\label{sec:utilization}

As we already mentioned before, the total glucose utilization $\U$ is the sum of two contributions~\cite[Sect. III-C-3]{dallaman07}: the insulin independent and the insulin dependent utilization, denoted by $\Uii$ and $\Uid$, respectively. More precisely, $\Uii$ represents the glucose uptake of the nervous system, which is thus assumed to be constant, i.e. $\Uii=\Fcns$. On the other hand, $\Uid$ depends on the glucose in tissues $\Gt$ as described by the Michaelis--Menten law
\begin{equation}
	\Uid(t) = \frac{\Vm(\X(t))\,\Gt(t)}{\Km(\X(t))+\Gt(t)},
\end{equation}
where $\Vm(x)=\Vmzero+x\,\Vmx$ and $\Km(x)=\Kmzero$ are linear functions of the insulin in the interstitial fluid $\X$ (pmol/L) governed by
\begin{equation}
	\dot{\X}(t) = -\pdueU\, \X(t) + \pdueU\,\big(\I(t)-\I_b\big).
\end{equation}
The coefficient $\pdueU$ (min$^{-1}$) is the rate of insulin action on the peripheral glucose utilization, $\I_b$ denotes the basal level of insulin concentration.

\rv{We finally point out that, at basal state, a direct relation holds between $\EGP$ and $\U$. Indeed, $\EGP_b=\U_b=\Uii+\Uid(0)$, see~\cite[Sect. III-C-3]{dallaman07}.}

\subsection{Insulin secretion}\label{sec:insulin_secretion}

The secretion of insulin by the pancreas is described by a set of four equations. We first observe that the insulin secretion is obtained by multiplying the portal insulin $\Ipo$ by the transfer rate $\gamma$ (min$^{-1}$) between portal vein and liver, i.e.
\begin{equation}
	\Se(t) = \gamma\,\Ipo(t).
\end{equation}
The portal insulin solves the following equation
\begin{equation}
	\Idot{po}=-\gamma\, \Ipo(t)+\Spo(t),
	\qquad\qquad\Ipo(0)=I_{\mathrm{po},b}
\end{equation}
where
\begin{equation}\label{eq:alpha}
	\begin{aligned}
		&\Spo(t) = 
		\begin{cases}
			\Y(t) + \K\, \dot{\G}(t) + \Se_b & \text{if }\dot{\G}>0\\
			\Y(t) + \Se_b & \text{if }\dot{\G}\le0,
		\end{cases}\\
		\qquad
		&\dot{\Y}(t) = 
		\begin{cases}
			-\alpha\,[\Y(t)-\beta\,(\G(t)-\h)]&\text{if } \beta\,(\G(t)-\h)\ge -\Se_b\\
			-\alpha\,\Y(t)-\alpha\, \Se_b&\text{if } \beta\,(\G(t)-\h)< -\Se_b
		\end{cases}\qquad \Y(0)=0,
	\end{aligned}
\end{equation}
and $\K$ (pmol/kg per mg/dL) is the pancreatic responsivity to the glucose rate of change, $\alpha$~(min$^{-1}$) denotes the delay between the glucose signal and the secretion of insulin by the pancreas, $\beta$ (pmol/kg/min per mg/dL) is the pancreatic responsivity to glucose. Finally, the constant $\h$~(mg/dL) is the threshold of glucose above which the $\beta-$cells begin the production of new insulin. Usually, $h=\G_b$. The symbol $\Se_b$ denotes the basal level of insulin secretion. 

\subsection{Glucose renal excretion}

The kidneys release glucose in plasma when its quantity decreases under a certain threshold $\kedue$ (mg/kg). This phenomenon is represented with a linear relation
\begin{equation}
	\E(t) = \begin{cases}
		\keuno(\Gp(t)-\kedue)
		&\text{if }\Gp(t)>\kedue\\
		0
		&\text{if }\Gp(t)\le\kedue.
	\end{cases}
\end{equation}
\rv{$\keuno$ (min$^{-1}$) is a rate constant describing renal glucose excretion. Since the model is formulated in terms of glucose masses (mg/kg), $\keuno$ represents a kinetic coefficient within this mass-based framework.}

\subsection{Delays}
    \rv{Similarly to~\cite{huard}, the considered model takes into account four physiological delays affecting the glucose-insulin dynamics. More precisely, the variable $\Id$, see~\eqref{eq:id}, is the delayed insulin signal, which directly affects the endogenous glucose production $\EGP$. The coefficient $\ki$, see~\eqref{eq:id}, denotes the delay between the insulin signal and the insulin action, while $\alpha$, see~\eqref{eq:alpha}, represents the delay between the glucose signal and the insulin secretion. Finally, the delay between the glucose ingestion and its assimilation is modeled by the rate of appearance $\RA$ of glucose into plasma, see~\eqref{eq:glucose}}.

\section{Parameters estimation}\label{sec:optimization}

In this section we describe the procedure we adopt to estimate part of the parameters involved in the system of ODEs. The estimation is performed by means of a constrained optimization algorithm on a dataset of sampled glucose curves.

Our attention is focused on the rate of appearance parameters, which describe the glucose kinetics in the gastrointestinal tract (see Section~\ref{sec:gastro_tract}) as well as $\Kempt$, and the endogenous glucose production parameters, appearing in the definition of $\EGP$, see~\eqref{eq:egp}. Such parameters have a direct effect on the blood glucose concentration $\G$ since they describe the absorption of glucose and its endogenous production. In Table~\ref{tab:intervals} we report the list of parameters we aim to estimate together with their lower and upper bounds, which have been extrapolated from previous studies (see e.g.~\cite{dallaman06,dallaman07}). We also report the parameters' initial value. All the other physical constants are fixed and listed in Table~\ref{tab:param_dallaman}~\cite[Tab. 1]{dallaman07}.  If a parameter is \red not reported in Table~\ref{tab:intervals} and Table~\ref{tab:param_dallaman}\black, then it is computed from some the other quantities. Regarding $\EGP$, we estimate its basal value $\EGP_b$ and then we compute $\kpuno$ by applying the relation in~\eqref{eq:egpb}.

\begin{table}\renewcommand{\arraystretch}{1.2}
	\centering
	\begin{tabular}{cccc}
		\textbf{Parameters} ($\ttheta$) & \textbf{Lower bound} & \textbf{Upper bound} & \textbf{Initial value} ($\ttheta_0$)\\
		\hline
		\multicolumn{4}{c}{\textit{Rate of Appearance parameters}}\\
		\hline
		$\Kmin$ & 1e-04 & 0.025 & 0.015\\
		$\Kmax$ & 0.035 & 0.1 & 0.0558\\
		$\Kabs$ & 0.01 & 0.3 & 0.057\\
		$\Kgri$ & 1e-05 & 0.1 & 0.049\\
		$\b$ & 0.65 & 0.995 & 0.85\\
		$\d$ & 1e-07 & 0.01 & 0.00018\\
		\hline
		\multicolumn{4}{c}{\textit{Endogenous Glucose Production parameters}}\\
		\hline
		$\EGP_b$ & 1.5 & 2.5 & 2\\
		$\kpdue$ & 1e-3 & 0.01 & 0.0021\\
		$\kptre$ & 1e-05 & 0.02 & 0.009\\
		$\kpquattro$ & 1e-04 & 0.1 & 0.0618\\
		$\ki$ & 1e-05 & 1e-03  & 0.0079\\
		\hline
	\end{tabular}
	\caption{Lower/upper bounds and initial value for parameters estimation. \rv{The bounds have been extrapolated from previous studies (see e.g.~\cite{dallaman06,dallaman07})}}
	\label{tab:intervals}
\end{table}

\begin{table}\renewcommand{\arraystretch}{1.2}
	\centering
	\begin{tabular}{cc||cc}
		\textbf{Parameter name} & \textbf{Value} & \textbf{Parameter} & \textbf{Value}\\
		\hline
		$\VG$ & 1.88 & $\Fcns$ & 1\\
		$\kuno$ & 0.065 & $\Vmzero $&2.50\\
		$\kdue$ & 0.079 & $\Vmx$ & 0.047\\
		$\VI$ & 0.05 & $\Kmzero$ & 225.59\\
		$\muno$ & 0.190 & $\pdueU$ & 0.0331\\
		$\mdue$ & 0.484 & $\K$ & 2.30\\
		$\mquattro$ & 0.194 & $\alpha$ & 0.050\\
		$\mcinque$ & 0.0304 & $\beta$ & 0.11\\
		$\msei$ & 0.6471 & $\gamma$ & 0.5\\
		$\HE_b$ & 0.6 & $\keuno$ & 1e-04\\
		$f$ & 0.90 & $\kedue$ & 339\\
		\hline
	\end{tabular}
	\caption{List of fixed parameters, \rv{see~\cite[Tab. 1]{dallaman07}}. During our simulation, we also fix the quantity of ingested glucose $D=50000$ mg and the body weight $\BW=78$ kg.}
	\label{tab:param_dallaman}
\end{table}

We denote by $(t_i,\mathsf{G}_i)$ the \textit{in vivo} sample points of a certain glucose concentration curve, measured at the time instant $t_i$ (min), for $i=1,\dots,N$, with $t_1=0$ corresponding to the basal state. Moreover $\ttheta=(\Kmin,\Kmax,\Kabs,\Kgri,\b,\d,\EGP_b,\kpdue,\kptre,\kpquattro,\ki)$ denotes the vector of parameters we are going to estimate and $\G(t;\ttheta)$ (mg/dl) is the associated glucose concentration curve resulting from the mathematical model. %We then represent by $\widetilde{\G}(t;\ttheta)$ the related curve measured in mmol/l, rescaled in such a way that $\widetilde{\G}(0;\ttheta)=0$.

In order to obtain the instance of $\ttheta$ providing a good fitting of the sample $(t_i,\mathsf{G}_i)$, we minimize the following loss function, where $\EGP(t;\ttheta)$ denotes the curve of $\EGP$ obtained from the mathematical model \red at the point $(t,\ttheta)$\black,
\begin{equation}\label{eq:loss}
	\loss(\ttheta) = \frac1N \sum_{i=1}^N \left[\mathsf{G}_i-{\G}(t_i;\ttheta)\right]^2 + \chi {\ind{z<0}}(\min\{\EGP(t;\ttheta)\}).
\end{equation}
More precisely, the first term is the mean squared error of the simulated glucose concentration at the sample points $t_i$ ($i=1,\dots,N$), while the second term forces $\EGP$ to be nonnegative by penalizing the error indicator by a constant $\chi=10^6$. Indeed, $\ind{z<0}$ is the characteristic function defined as
\begin{equation}
	\ind{z<0}(x) = \begin{cases}
		1 & \text{if }x<0\\
		0 & \text{otherwise}.
	\end{cases}
\end{equation}
\rv{In other words, whenever $\EGP$ assumes negative values during the minimization process, we have that
$$
{\ind{z<0}}(\min\{\EGP\})=1,
\qquad\text{and}\qquad
\loss(\ttheta) = \frac1N \sum_{i=1}^N \left[\mathsf{G}_i-{\G}(t_i;\ttheta)\right]^2 + \chi,
$$
with $\chi$ dominating over the mean squared error. This penalization avoids the choice of parameters values giving a good fitting of the sample points, but an inadmissible profile of $\EGP$.}

\rv{As the parameters are estimated within the ranges reported in Table~\ref{tab:intervals} (\rv{second and third column}), }the minimization process is carried out by means of the constrained optimization function \texttt{fmincon} provided by MATLAB, which implements the interior point algorithm discussed in~\cite{nocedal99,nocedal06,nocedal14}. The process is initialized by setting $\ttheta=\ttheta_0$ with the values reported in Table~\ref{tab:intervals} (\rv{fourth column}) and terminates when $\loss(\ttheta)<\tol$ or $\ttheta$ varies less than a prescribed tolerance $\tol$. In particular, we fix $\tol=$1e-10. As additional stopping criteria we adopt the maximum number of iterations and the maximum number of evaluations of $\loss$, which must be smaller than~$500$ times the size of $\ttheta$ (i.e. 5500 iterations). At each iteration of \texttt{fmincon}, the ODE system is solved on a time grid with step $\dt=0.05$.

\section{Results}\label{sec:results}

\subsection{Our dataset}
The database consists of 35 glycemic response curves obtained by 35 healthy subjects (range of age: 20-40) provided by the Laboratory of Neurobiology and Integrated Physiology (University of Pavia, Italy).

Each volunteer received precise behavioral instructions to avoid distorting the test results. In particular, they were advised not to consume alcohol in the previous 12 hours, smoke cigarettes, or participate in strenuous physical exercise in the previous days. Each subject received 50 g of glucose orally (Glucose Sclavo, diagnostic 75 g/150 ml) after fasting for at least 12 hours, according to the procedure defined by the International Organization for Standardization \cite{iso}. Blood glucose was monitored using the Lifescan One Touch UltraEasy\textsuperscript{\textcopyright} system, which uses a glucose oxidase biosensor as a dosing method \cite{Brouns2005,Wolever2003}. Specifically, blood glucose was measured before glucose ingestion (corresponding to basal or fasting glycemia) and then every 15 minutes after ingestion up to 2 hours later. %\red At the end of each measurement, the blood glucose values obtained, expressed in mg/dl, were converted to mM/L, and then the basal blood glucose value was subtracted from each value.\black
%
%\
%
%\NOTE{Questa ultima frase forse \`e da togliere: noi facciamo i conti su curve non scalate e non normalizzate.}

\subsection{Estimated curves}

\begin{figure}
	\centering
	\centering
	\subfloat[]{\includegraphics[width=0.4\linewidth]{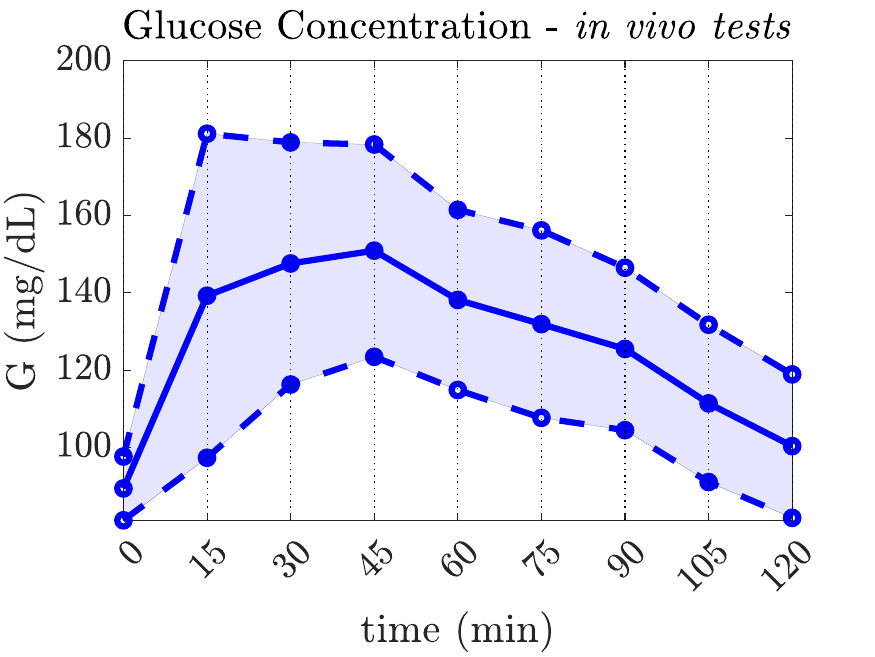}}\qquad
	\subfloat[]{\includegraphics[width=0.4\linewidth]{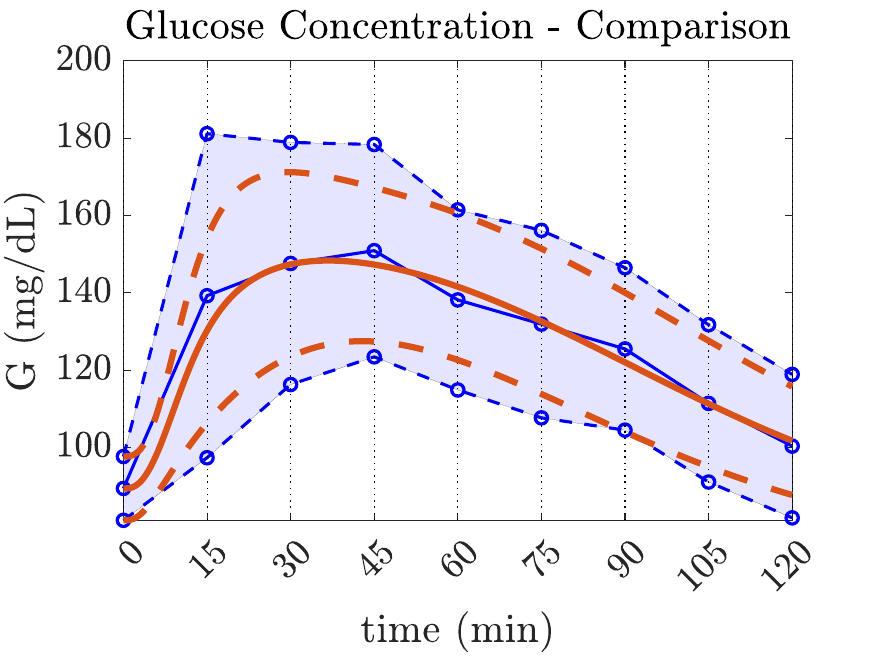}}
	\caption{Comparison between \textit{in vivo} and estimated glucose concentration curves over the 35 considered subjects. \textsc{(a)} Mean and deviation of the \textit{in vivo} glucose curves of our database. \textsc{(b)} Comparison between mean and deviation of \textit{in vivo} curves (blue lines) and estimated curves (orange lines).}
	\label{fig:sample_vs_model}
\end{figure}

The model simulation allows for the description of the physiological events that occur in the postprandial state following the ingestion of glucose solution by the healthy subjects in the database. As previously described, 35 healthy subjects consumed a solution containing 50 g of glucose within 15 minutes, after which glucose concentration levels were measured every 15 minutes for 2 hours. Figure~\ref{fig:sample_vs_model} shows a comparison between the data measured in vivo and the estimates obtained using the model. The model accurately estimates the glucose curves measured in vivo. Moreover, Table~\ref{tab:params} reports the average value and the standard deviation of each estimated parameter.

\begin{table}\renewcommand{\arraystretch}{1.2}
	\centering
	\begin{tabular}{cc}
		\textbf{Parameter name} & \textbf{Mean} $\pm$ \textbf{SD}\\
		\hline
		\multicolumn{2}{c}{\textit{Rate of Appearance parameters}}\\
		\hline
		$\Kmin$ & 0.0086 $\pm$ 0.0102\\
		$\Kmax$ & 0.0831 $\pm$ 0.0246\\
		$\Kabs$ & 0.2266 $\pm$ 0.0943\\
		$\Kgri$ & 0.0785 $\pm$ 0.0262\\
		$\b$ &  0.7391 $\pm$ 0.0976\\
		$\d$ & 0.0050 $\pm$ 8.1557e-06\\
		\hline
		\multicolumn{2}{c}{\textit{Endogenous Glucose Production parameters}}\\
		\hline
		$\EGP_b$ & 1.8871 $\pm$ 0.4055\\
		$\kpdue$ & 0.0044 $\pm$ 0.0040\\
		$\kptre$ & 0.0143 $\pm$ 0.0072\\
		$\kpquattro$ & 0.0251 $\pm$ 0.0376\\
		$\ki$ & 6.7815e-04 $\pm$ 4.1366e-04\\
		\hline
	\end{tabular}
    \caption{\rv{Mean value and standard deviation for the parameters we estimated starting from our dataset.}}
	%\caption{Mean value and standard deviation for each parameter estimated on our dataset.}
	\label{tab:params}
\end{table}

\begin{figure}
	\centering
	\subfloat[]{\label{fig:stats_a}\includegraphics[width=0.35\linewidth]{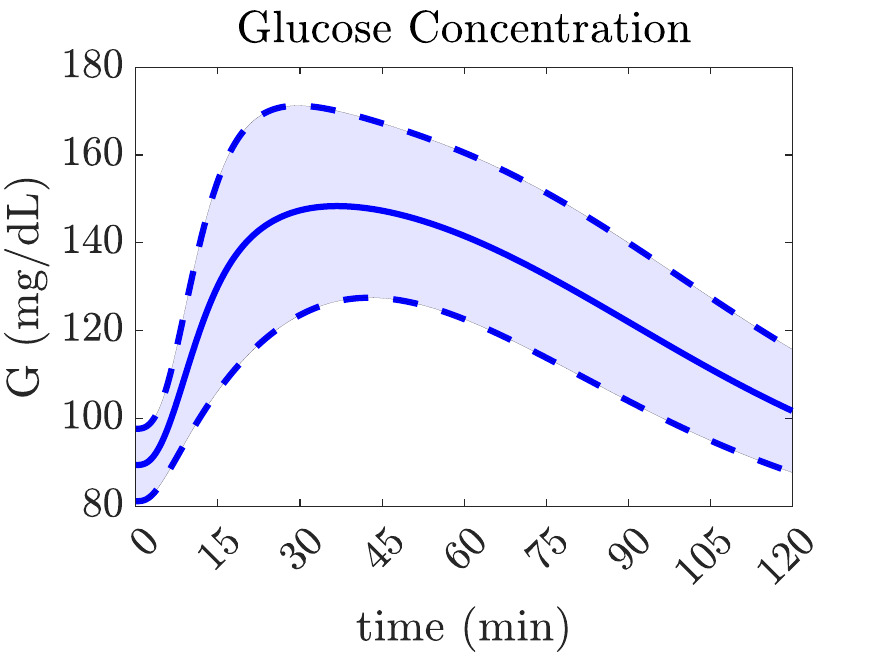}}
	\subfloat[]{\label{fig:stats_b}\includegraphics[width=0.35\linewidth]{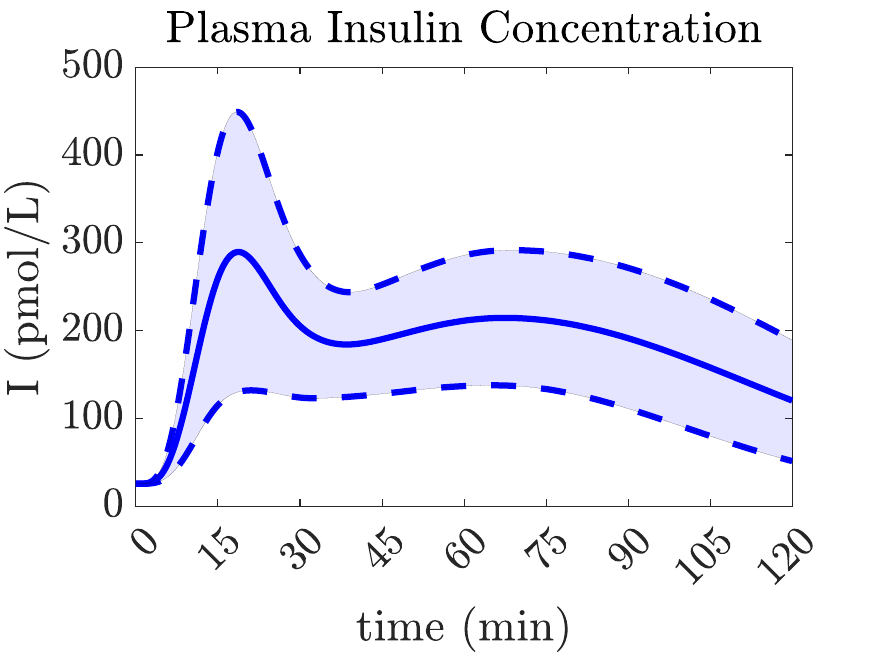}}
	\subfloat[]{\label{fig:stats_c}\includegraphics[width=0.35\linewidth]{stat_egp.pdf}}\\
	\subfloat[]{\label{fig:stats_d}\includegraphics[width=0.35\linewidth]{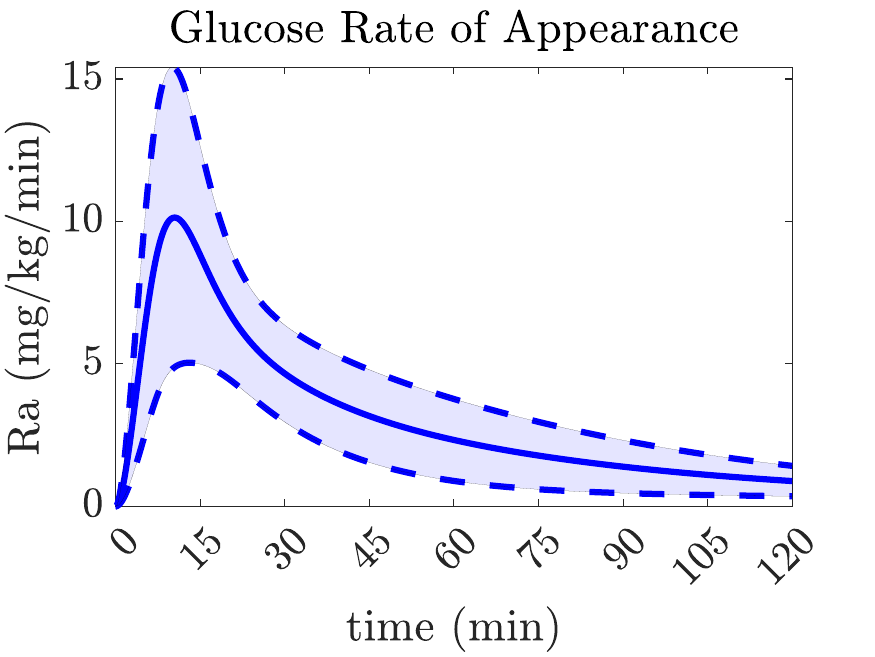}}
	\subfloat[]{\label{fig:stats_e}\includegraphics[width=0.35\linewidth]{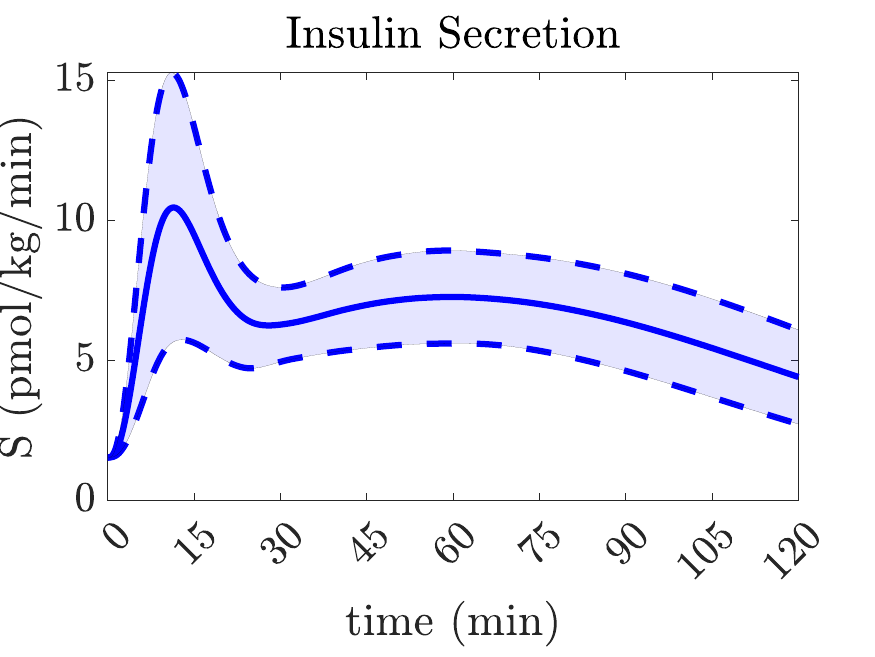}}
	\subfloat[]{\label{fig:stats_f}\includegraphics[width=0.35\linewidth]{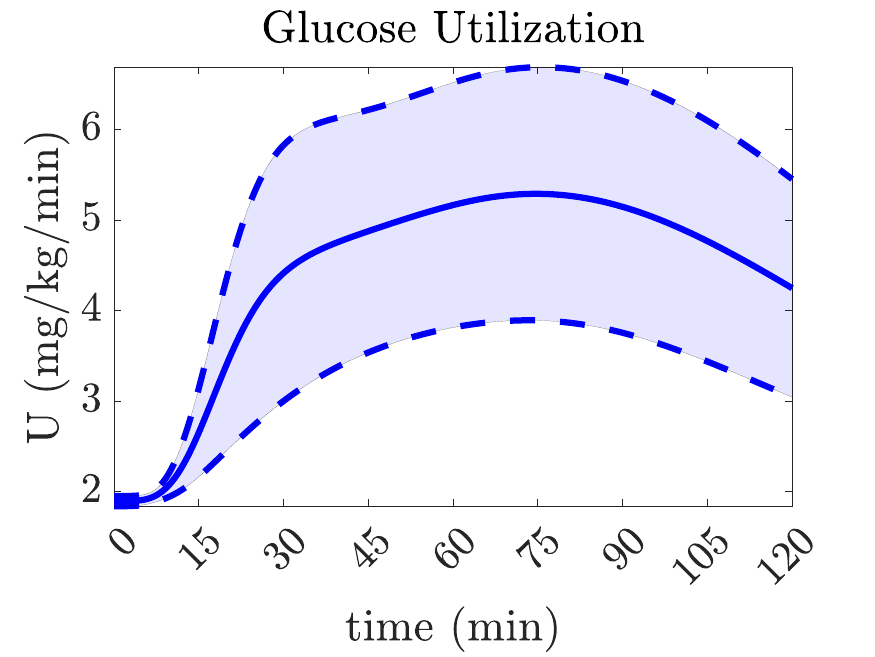}}
	\caption{Quantities of \textsc{(a)} glucose concentration, \textsc{(b)} plasma insulin concentration, \textsc{(c)} endogenous glucose production, \textsc{(d)} glucose rate of appearance, \textsc{(e)} insulin secretion, and \textsc{(f)} glucose utilization. The solid line represents the mean of 35 subjects, and the shaded area represents the SD.}
	\label{fig:stats}
\end{figure}

The dynamic interaction between glucose absorption, endogenous production, insulin secretion, and glucose utilization during the postprandial state is plotted in Figure~\ref{fig:stats}. More precisely, plasma insulin concentration (Figure~\ref{fig:stats_b}) exhibits a pattern similar to that of glucose (Figure~\ref{fig:stats_a}), with a delayed peak compared to glucose. This delay is due to the time required for insulin to reach the plasma circulation after being released by pancreatic $\beta-$cells in response to the rise in glucose levels. Notably, insulin secretion (Figure~\ref{fig:stats_e}) shows an early and large peak that precedes the plasma one and then progressively reduces as glucose concentration decreases. 
Similarly, endogenous glucose production (Figure~\ref{fig:stats_c}) undergoes a rapid reduction during the initial minutes, driven by the effect of insulin, followed by a gradual recovery, \red as the trends \black of glucose and insulin in the bloodstream.

The profile of glucose utilization (Figure~\ref{fig:stats_f}) is similar to that of insulin, with a progressive increase up to a plateau. The variability reflects the differences between individuals in insulin sensitivity and cellular glucose metabolism. The rate of glucose appearance $\RA$ (Figure~\ref{fig:stats_d}) decreases after an initial peak: this behavior is consistent with the rapid absorption of glucose from the intestine, which immediately follows a meal. This trend reflects the systemic clearance of glucose from the plasma.

The shaded regions, common across all panels, represent the standard deviation (SD) emphasizing the \red{model's} ability to capture not only the mean physiological processes but also the inter-individual variability. These findings validate the model's capacity to accurately simulate the key physiological processes involved in glucose metabolism, with particular emphasis on the temporal dynamics and individual variability of the responses.

\subsection{Three behaviors}

In this section, we describe the individual glucose tolerance of the subjects, along with its relation with the dynamics of the entire system. %and we analyze its relation with the estimated parameters.

\begin{figure}
	\centering
	\subfloat[\label{fig:three_behaviors_g}]{\includegraphics[width=0.35\linewidth]{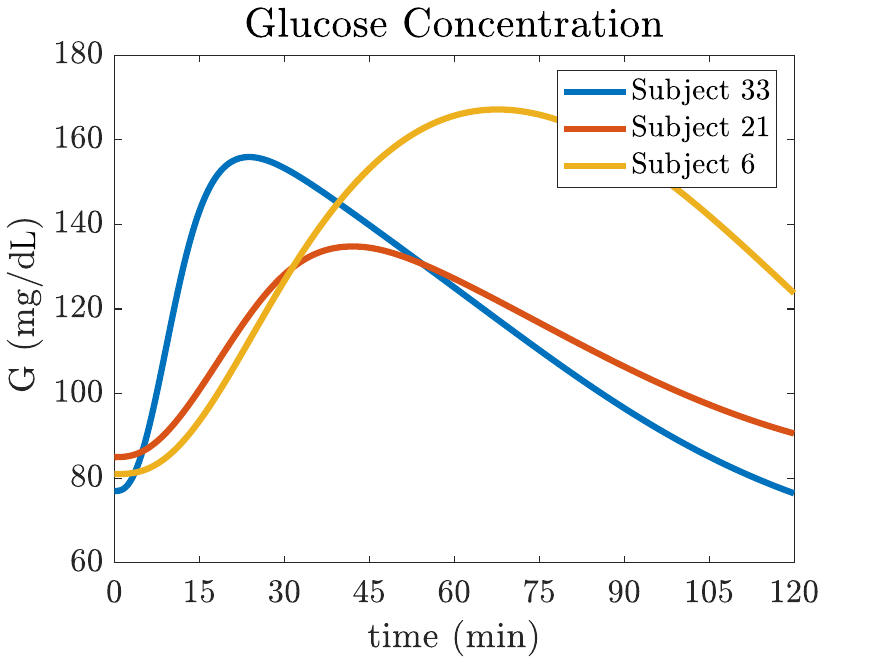}}
	\subfloat[\label{fig:three_behaviors_i}]{\includegraphics[width=0.35\linewidth]{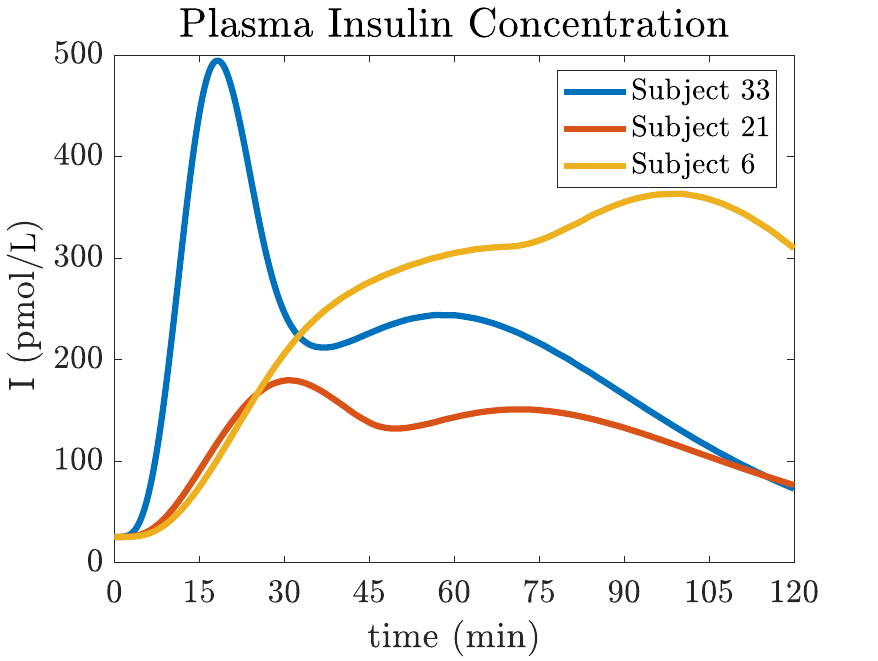}}
	\subfloat[\label{fig:three_behaviors_egp}]{\includegraphics[width=0.35\linewidth]{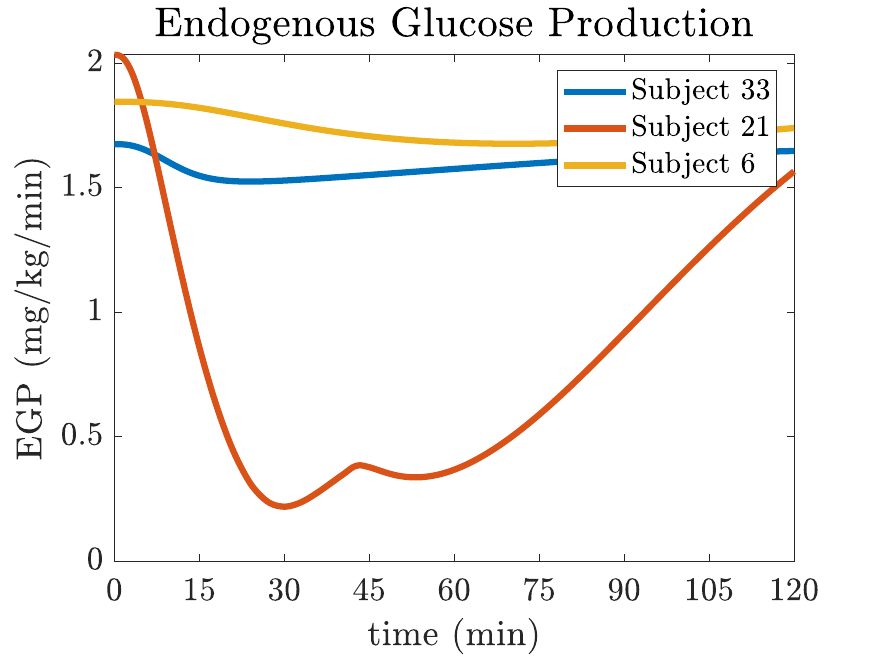}}\\
	\subfloat[\label{fig:three_behaviors_ra}]{\includegraphics[width=0.35\linewidth]{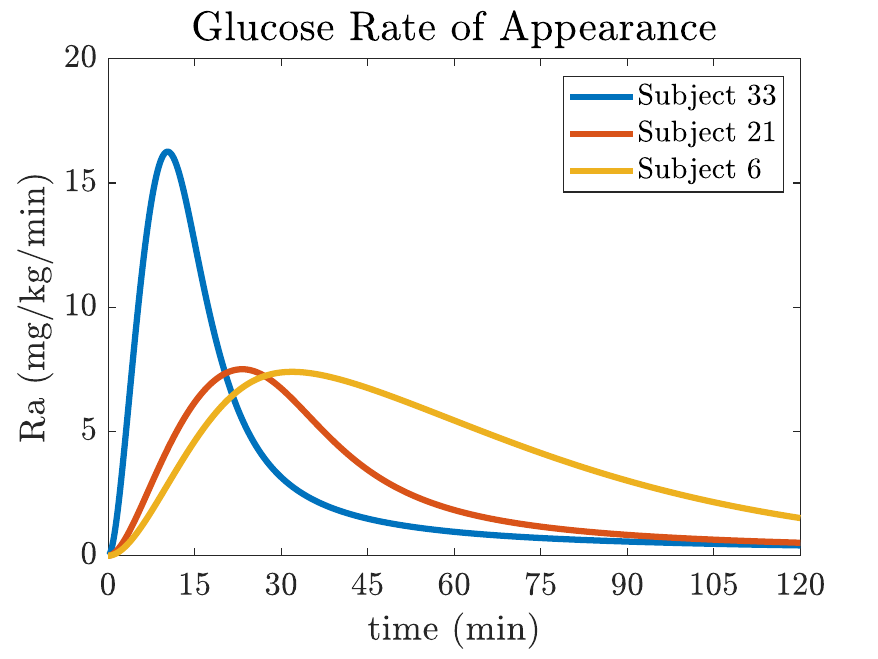}}
	\subfloat[\label{fig:three_behaviors_s}]{\includegraphics[width=0.35\linewidth]{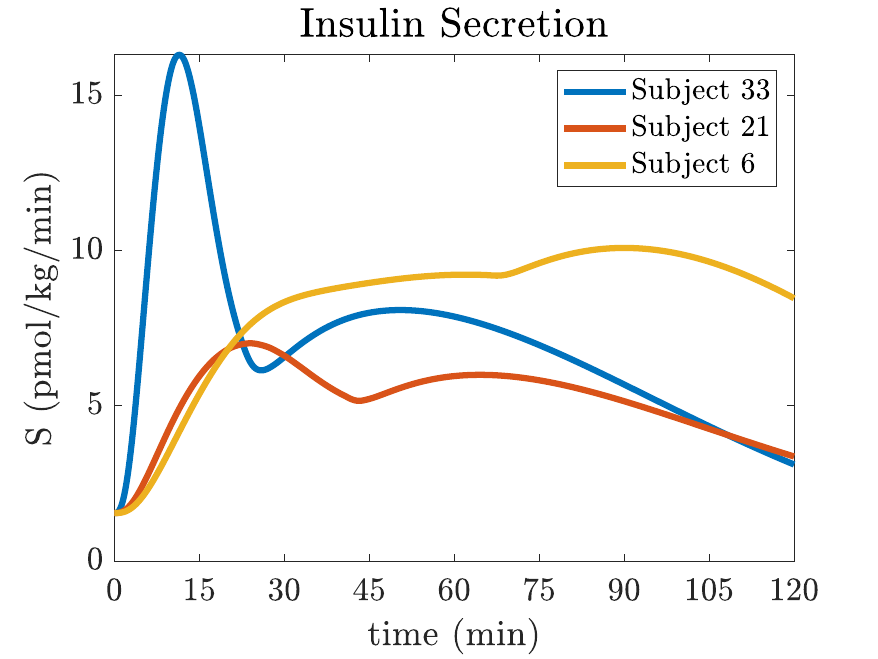}}
	\subfloat[\label{fig:three_behaviors_u}]{\includegraphics[width=0.35\linewidth]{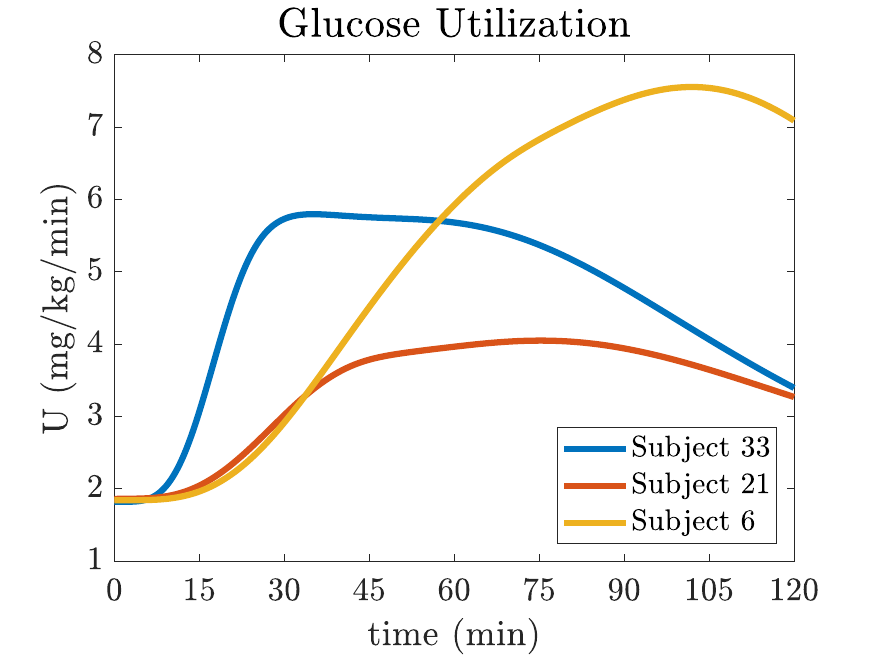}}
	\caption{Quantities of \textsc{(a)} glucose concentration, \textsc{(b)} plasma insulin concentration, \textsc{(c)} endogenous glucose production, \textsc{(d)} glucose rate of appearance, \textsc{(e)} insulin secretion, and \textsc{(f)} glucose utilization of subjects 33 (blue line), 21 (red line), 6 (yellow line).}
	\label{fig:three_behaviors}
\end{figure}

More precisely, we identified three distinct behaviors as represented in Figure~\ref{fig:three_behaviors} for three reference subjects. From Figure~\ref{fig:three_behaviors_g}, we observe that such classification might depend on the time of glycemic peak. In fact, the glycemic peak of Subject 33 (blue line) occurred within 30 minutes after ingestion, indicating rapid glucose absorption, as also shown by the associated $\RA$ and $\U$ curves (see Figures~\ref{fig:three_behaviors_ra} and~\ref{fig:three_behaviors_u}, respectively). However, looking at Subject 21 (red line), the glycemic peak occurred between 30 and 50 minutes after ingestion. This behavior indicats a moderate glycemic response, as seen in the appearance rate and the glucose utilization plots (see Figures~\ref{fig:three_behaviors_ra} and~\ref{fig:three_behaviors_u}, respectively). Subject 6 (yellow line) showed a completely different behavior, as the glycemic peak was more delayed, occurring 60 minutes after ingestion. The related glucose rate of appearance (Figure~\ref{fig:three_behaviors_ra}) exhibited a gradual and slow glucose uptake correlated with a low ability to utilize ingested glucose, see Figure~\ref{fig:three_behaviors_u}.

As represented in Figures~\ref{fig:three_behaviors_i} and~\ref{fig:three_behaviors_s}, the glucose dynamics affected the insulin response. Subject 33 had a high pancreatic response, resulting in a significant increase in insulin concentration to compensate for the immediate glucose peak. The immediate action of insulin led to a decrease in endogenous glucose production (Figure~\ref{fig:three_behaviors_egp}). Subject 21 showed regular and moderate insulin secretion, which is correlated with a gradual suppression of gluconeogenesis. Due to its delayed glucose peak, Subject 6 exhibited a slow and progressive increase of insulin secretion and action, resulting in a remarkable contribution in reducing the endogenous glucose production.

The three subjects displayed different metabolic reactions, reflecting the biological variability of glycemic regulation.

\subsubsection{Analysis of rate of appearance parameters}

\begin{figure}
	\centering
	\subfloat[]{\includegraphics[width=0.4\linewidth]{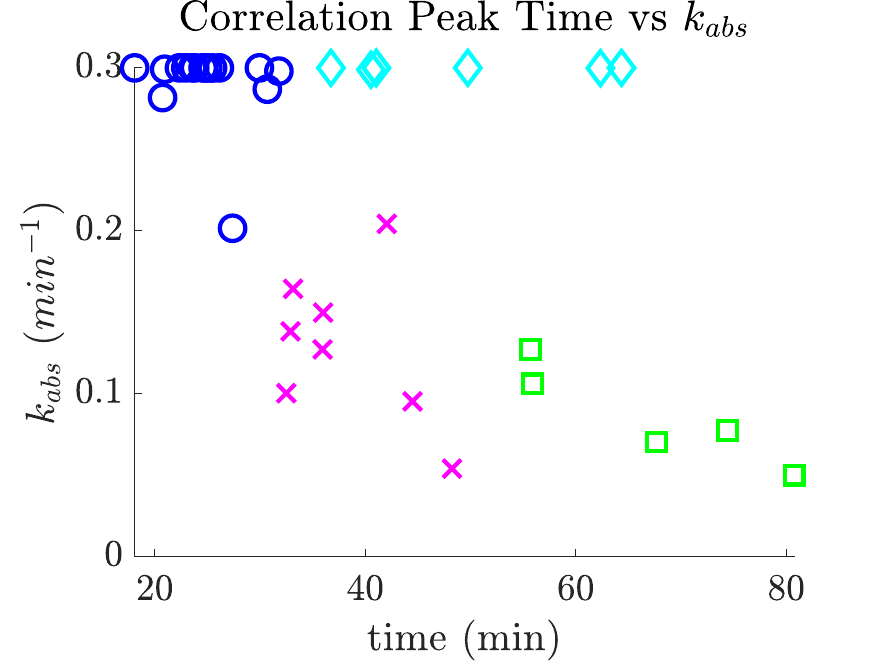}}\qquad
	\subfloat[]{\includegraphics[width=0.19\linewidth]{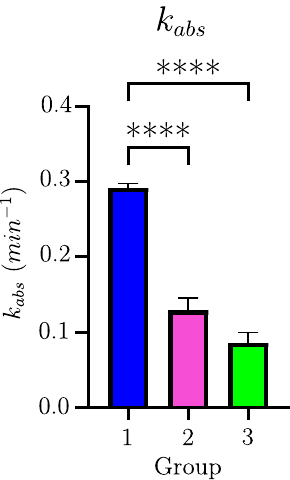}}
	\caption{Correlation plots between time of glycemic peak and $\Kabs$ (rate of absorption). \textsc{(a)} Scatter plot of $\Kabs$ for Group 1 (blue circles), Group 2 (magenta crosses), Group 3 (green squares), and six outliers (cyan diamond). \textsc{(b)} Statistical comparison of $\Kabs$ among the 3 groups. Statistical significance (One--Way ANOVA followed by Bonferroni \emph{post-hoc} test): $p < 0.05$ (*); $p < 0.01$ (**); $p < 0.001$ (***); $p < 0.0001$ (****).}
	\label{fig:corr_kabs}
\end{figure}

\begin{figure}
	\centering
	\subfloat[]{\includegraphics[width=0.4\linewidth]{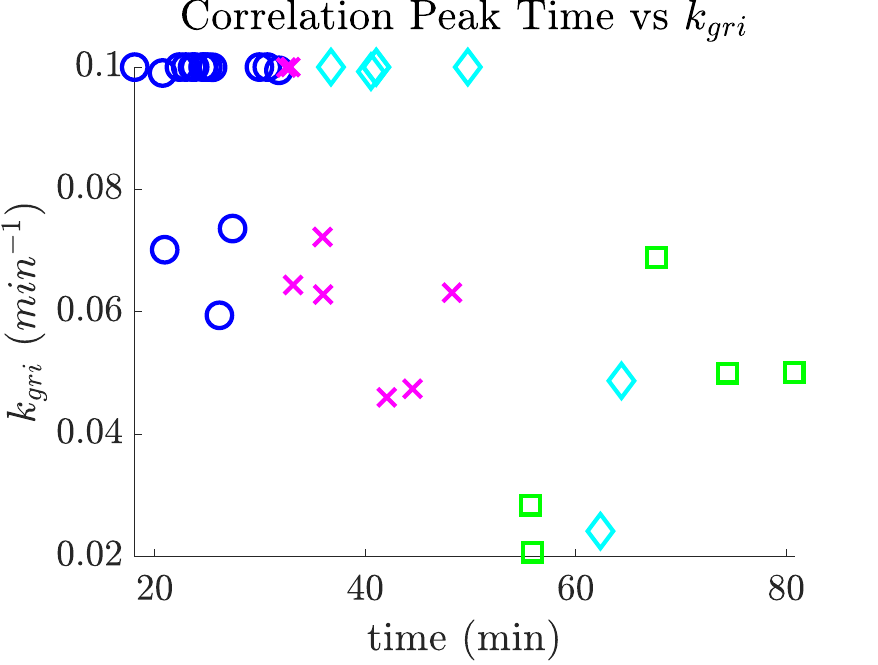}}\qquad
	\subfloat[]{\includegraphics[width=0.2\linewidth]{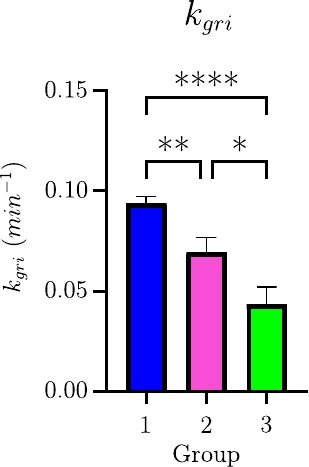}}
	\caption{Correlation plots between time of glycemic peak and $\Kgri$ (rate of grinding). \textsc{(a)} Scatter plot of $\Kgri$ for Group 1 (blue circles), Group 2 (magenta crosses), Group 3 (green squares), and six outliers (cyan diamond). \textsc{(b)} Statistical comparison of $\Kgri$ among the 3 groups. Statistical significance (One--Way ANOVA followed by Bonferroni \emph{post-hoc} test): $p < 0.05$ (*); $p < 0.01$ (**); $p < 0.001$ (***); $p < 0.0001$ (****).}
	\label{fig:corr_kgri}
\end{figure}

Following the previous discussion, we analyze the correlation between the three behaviors and the estimated rate of appearance parameters. To this end, we subdivide the database into three groups based on the time of the glucose peak. Group 1 corresponds to the subject with an early peak (30 minutes after glucose ingestion); Group 2 corresponds to a medium peak (between 30 and 50 minutes after ingestion); and Group 3 corresponds to a late peak (after 50 minutes). 
The analysis aims at exploring the relationship between peak time and two key kinetic parameters: the absorption rate constant ($\Kabs$) and the grinding rate constant ($\Kgri$) (see Table~\ref{tab:groups}, and Figures~\ref{fig:corr_kabs} and ~\ref{fig:corr_kgri}). We have six outliers, which will be subsequently analyzed in detail.

\begin{table}\renewcommand{\arraystretch}{1.2}
	\centering
	\begin{tabular}{ccccc}
		& \# Subjects & Peak Time & $\Kabs$ & $\Kgri$\\
		\hline
		Group 1 & 16 & $24.8938 \pm 0.9295$ & $0.2917 \pm 0.0062$ & $0.0939 \pm 0.0033$\\
		Group 2 & 8 & $38.1625 \pm 2.1192$ & $0.1291 \pm 0.0164$ & $0.0695 \pm 0.0073$\\
		Group 3 & 5 & $66.8800 \pm 4.9753$ & $0.0861 \pm 0.0136$ & $0.0436 \pm 0.0086$\\
		%Outliers & 6 & -- & -- & --\\
		%\hline
		%Total & 35 & -- & -- & --\\
		\hline
	\end{tabular}
	\caption{Mean values and standard deviation of peak time, $\Kabs$, and $\Kgri$ for subjects divided into three groups based on the time of the glycemic peak: Group~1 (peak time $<30$ min), Group~2 (peak time between $30-50$ min), and Group~3 (peak time $>50$ min).}
	\label{tab:groups}
\end{table}

% SCATTER PLOTS PER VALORI STIMATI con OUTLIERS

\begin{figure}
	\centering
	\subfloat[]{\includegraphics[width=0.4\linewidth]{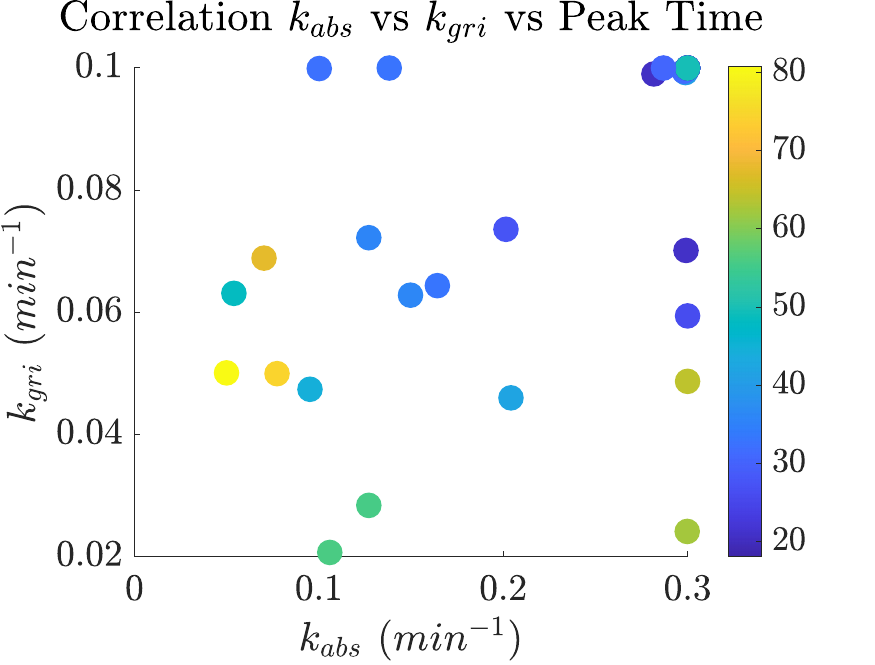}}\qquad
	\subfloat[]{\includegraphics[width=0.4\linewidth]{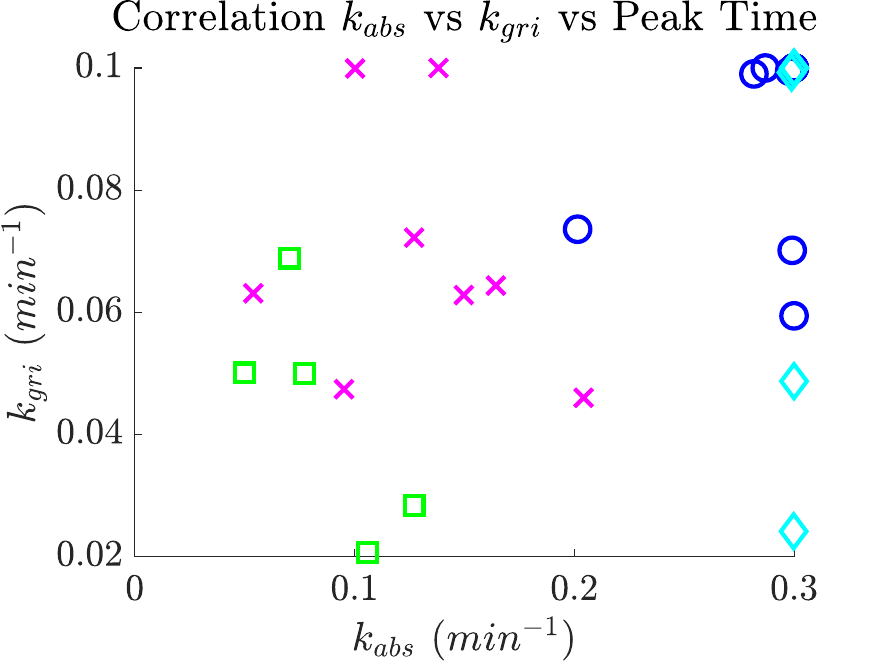}}
	\caption{Correlation plots between $\Kabs$, $\Kgri$ and the time of glycemic peak. \textsc{(a)}~Scatter plot of $\Kabs$ and $\Kgri$, where the color scale represents the peak time (min). \textsc{(b)}~Scatter plot of $\Kabs$ and $\Kgri$ for Group 1 (blue circles), Group 2 (magenta crosses), Group 3 (green squares), and six outliers (cyan diamond).}
	\label{fig:correlations}
\end{figure}

In particular, the subjects of Group 1 shows the largest absorption rate ($\Kabs = 0.2917 \pm 0.0062$ min$^{-1}$) and the largest grinding rate ($\Kgri = 0.0939 \pm 0.0033$ min$^{-1}$), indicating a rapid metabolic response. For subjects in Group 2, the absorption rate is markedly lower than in Group 1 ($\Kabs = 0.1291 \pm 0.0164$ min$^{-1}$), and their grinding rate is also reduced ($\Kgri = 0.0695 \pm 0.0073$ min$^{-1}$), suggesting slower glucose uptake and metabolic processing.
Instead, subjects in Group~3 exhibit the lowest values for both kinetic parameters ($\Kabs = 0.0861 \pm 0.0136$ min$^{-1}$, $\Kgri = 0.0436 \pm 0.0086$ min$^{-1}$), indicating a significantly slower absorption and grinding process compared to the other groups.
These findings suggest a clear inverse relationship between the time of glycemic peak and the kinetic parameters governing glucose absorption and processing. Individuals who exhibited an earlier glycemic peak tend to have higher metabolic efficiency, characterized by faster glucose absorption (larger $\Kabs$) and increased glucose processing (larger $\Kgri$). However, individuals with delayed glycemic peaks show a markedly slower metabolic response, indicating physiological variations in glucose handling that may be essential to evaluate metabolic health and disease risk.
Furthermore, the correlations between $\Kabs$, $\Kgri$, and the glycemic peak time (see  Figure~\ref{fig:correlations}) confirm an inverse relationship between peak time and metabolic efficiency, highlighting glucose absorption and processing variability.

% SPOSTATO QUESTA FRASE SOTTO, quando si parla effettivamente di outliers. Poi vedere se non \`e il caso di riscrivere una frasetta anche qui
%In addition, six subjects were identified as outliers, and their data were excluded from the grouped statistical analyses for further investigation.
%
%\NOTE{Bisogna parlare degli outlier, magari qui bisogna aggiungere una frase di connessione.}

\begin{figure}
	\centering
	\subfloat[]{\includegraphics[width=0.4\linewidth]{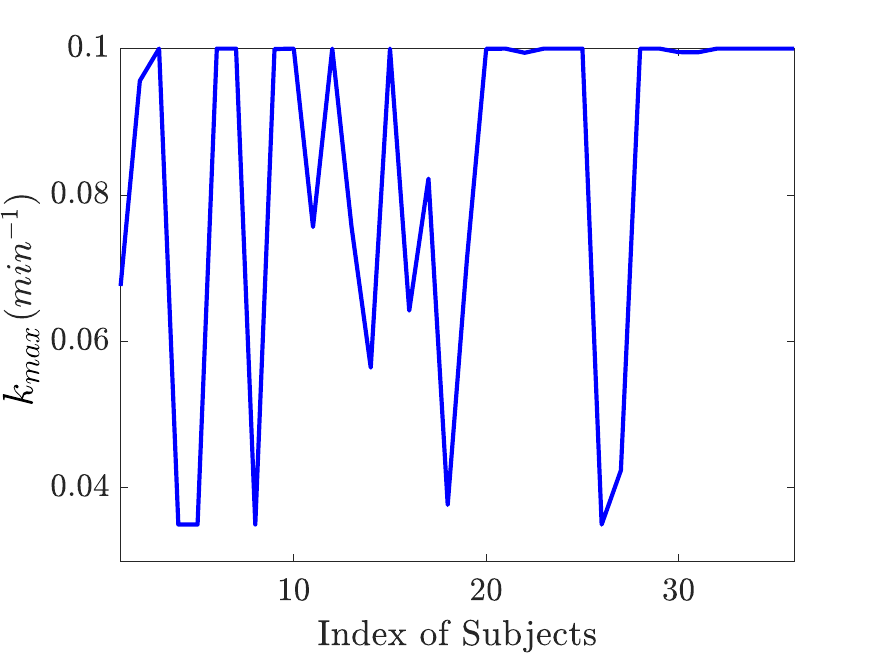}}\qquad
	\subfloat[]{\includegraphics[width=0.4\linewidth]{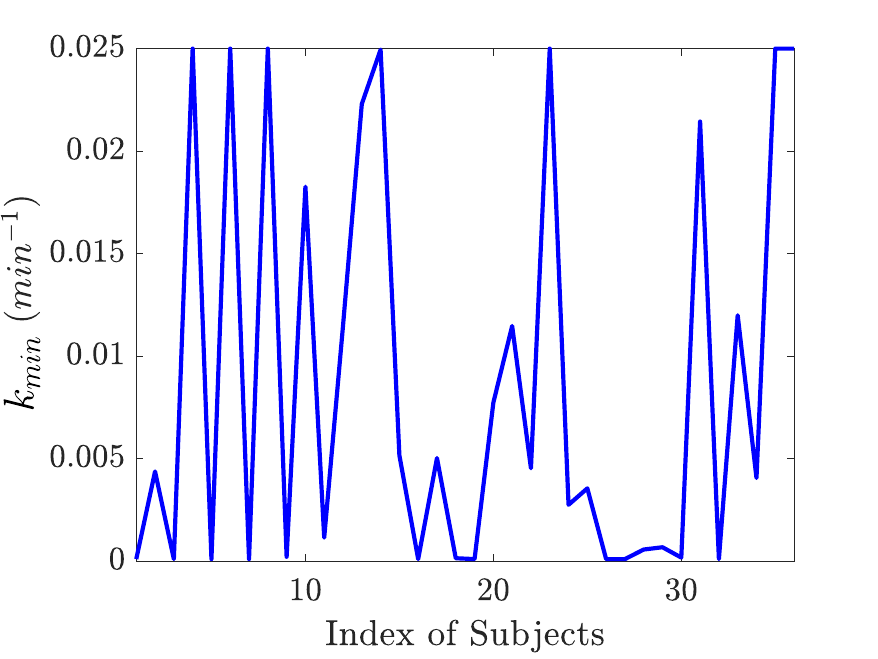}}\\
	
	\
	
	\subfloat[]{\includegraphics[width=0.4\linewidth]{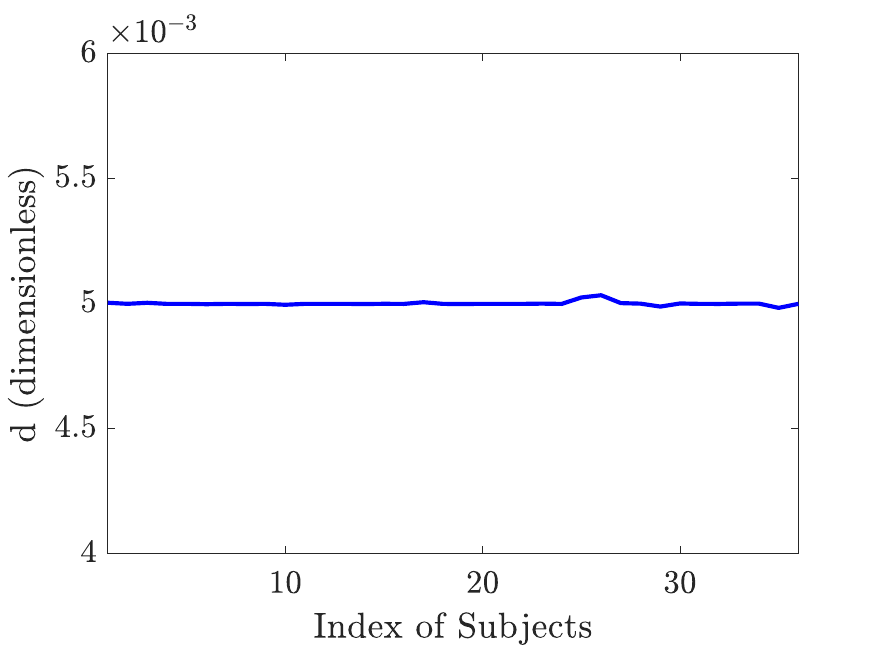}}\qquad
	\subfloat[]{\includegraphics[width=0.4\linewidth]{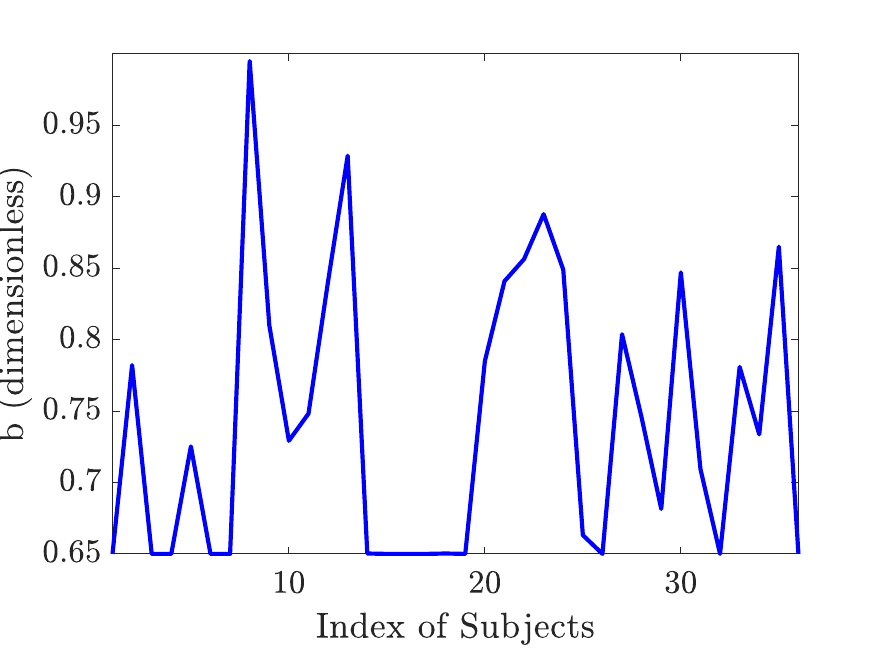}}
	\caption{Gastric emptying parameters with respect to the index of subjects. \textsc{(a)}~$\Kmax$, \textsc{(b)}~$\Kmin$, \textsc{(c)}~$\d$, \textsc{(d)}~$\b$.}
	\label{fig:empt}
\end{figure}

We analyzed also other rate of appearance parameters (see Figure~\ref{fig:empt}). $\Kmax$, $\Kmin$, $\d$, and $\b$ describe different aspects of the gastric emptying process, capturing its variability and dynamics across subjects. More precisely, $\Kmax$ and $\Kmin$ represent the maximum and minimum rate of gastric emptying, respectively, and reflect normal biological variability. The parameter $\b$, which represents the percentage amount of glucose remaining in the stomach when the emptying rate has decreased to halfway between $\Kmax$ and $\Kmin$ (see~\eqref{eq:kempt} and Figure~\ref{fig:kempt}), appears to be correlated with both the timing and maximum value of glucose concentration.. 

\begin{figure}
	\centering
	\subfloat[]{\includegraphics[width=0.38\linewidth]{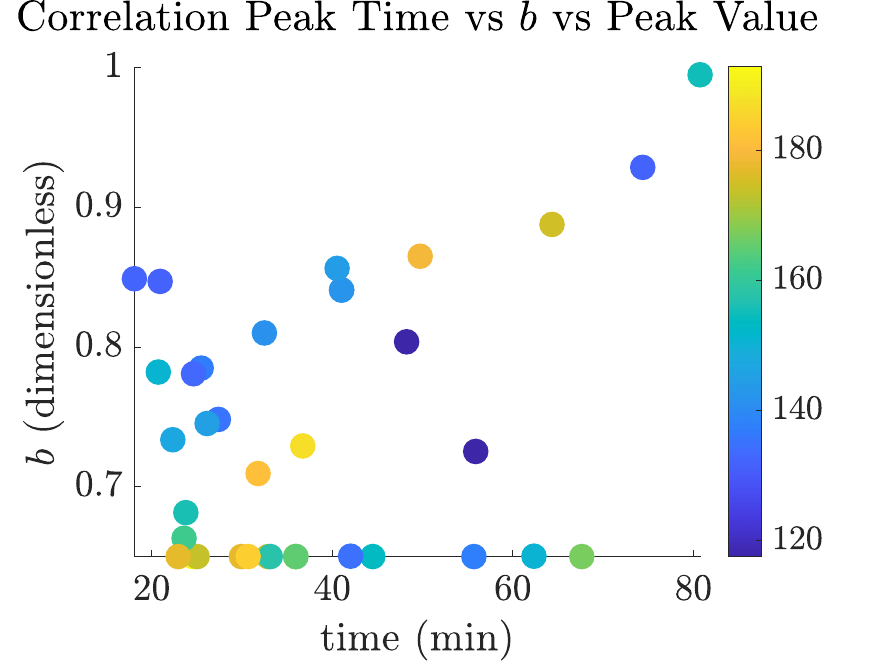}}
	\subfloat[]{\includegraphics[width=0.385\linewidth]{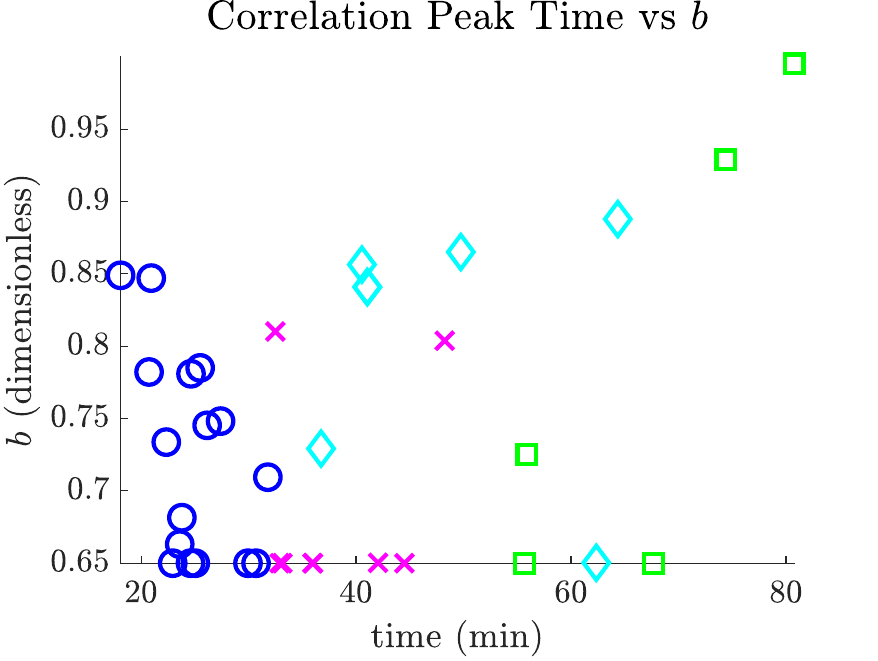}}
	\subfloat[]{\includegraphics[width=0.25\linewidth]{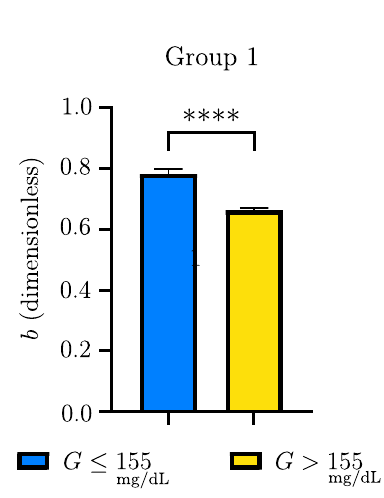}}
	\caption{Correlation plots between the time of glycemic peak, $b$ and the glycemic peak value. \textsc{(a)}~Scatter plot of Peak Time, $\b$, and Peak Value, where the color scale represents the peak value ($mg/dl$). \textsc{(b)}~Scatter plot of Peak Time and $\b$ for Group 1 (blue circles), Group 2 (magenta crosses), Group 3 (green squares), and six outliers (cyan diamond). \textsc{(c)}~Comparing the parameter $\b$ Group 1 for subjects with glucose values $\G\le155$ mg/dL (blue) and $\G>155$ mg/dL (yellow). Statistical significance (One--Way ANOVA followed by Bonferroni \emph{post-hoc} test): $p < 0.05$ (*); $p < 0.01$ (**); $p < 0.001$ (***); $p < 0.0001$ (****).}
	\label{fig:b}
\end{figure}

When $\b$ assumes a large value, the emptying rate $\Kempt$ decreases faster to $\Kmin$, meaning that the emptying process is slow. This implies that the glucose assimilation is prolonged over time. On the other hand, a small value of $\b$, results in a slower decrease of $\Kempt$ to $\Kmin$ so that the emptying process and glucose assimilation are faster. These observations are evident for subjects belonging to Group~1. Indeed, the subjects with $\G\le155$ mg/dL show a statistically larger $\b$ compared to subjects with $\G>155$ mg/dL, see Figure~\ref{fig:b}. Regarding Group~2 and Group~3, no clear classification could be identified in dependence on gastric emptying. 

The parameter $\d$, characterizing the recovery of $\Kempt$ to $\Kmax$, remains constant across the group of subjects, indicating a stable parameter in the gastric emptying process.

In any case, the correlation between $\b$ and the maximum of glucose concentration observed in Group~1 further supports the idea that gastric emptying dynamics plays a crucial role in shaping postprandial glycemic responses.

\begin{figure}
	\centering
	\subfloat[]{\includegraphics[width=0.348\linewidth]{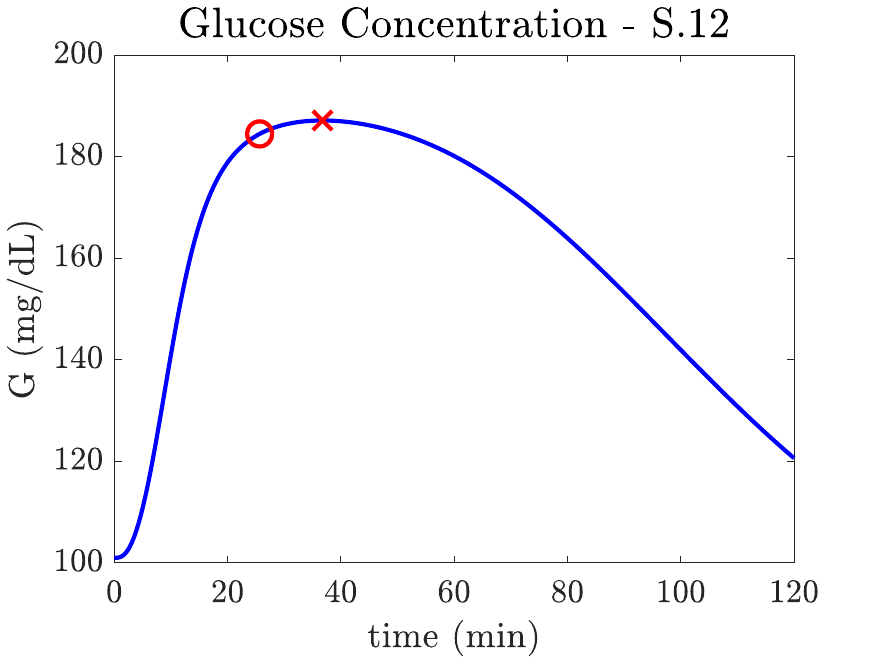}}
	\subfloat[]{\includegraphics[width=0.35\linewidth]{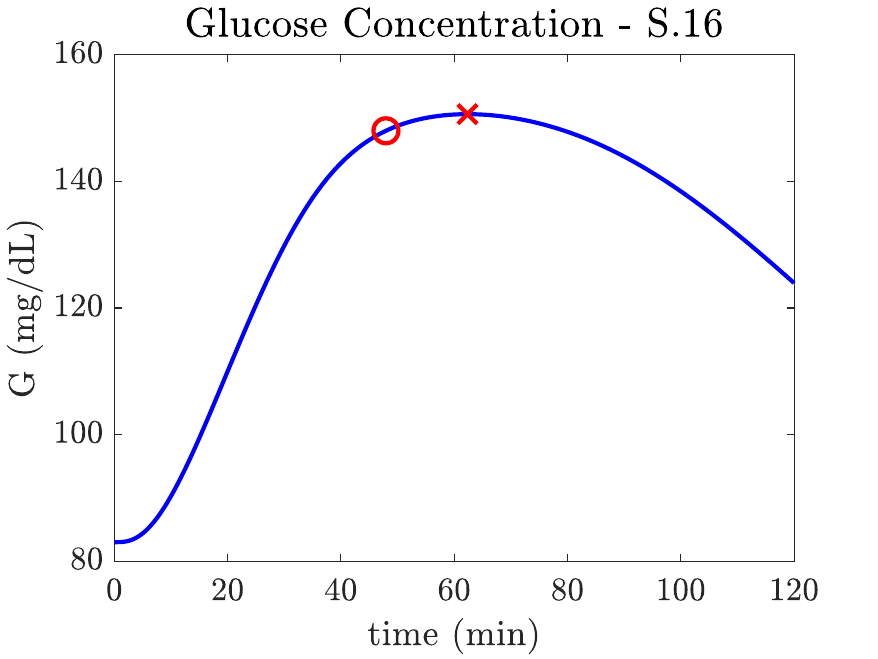}}
	\subfloat[]{\includegraphics[width=0.35\linewidth]{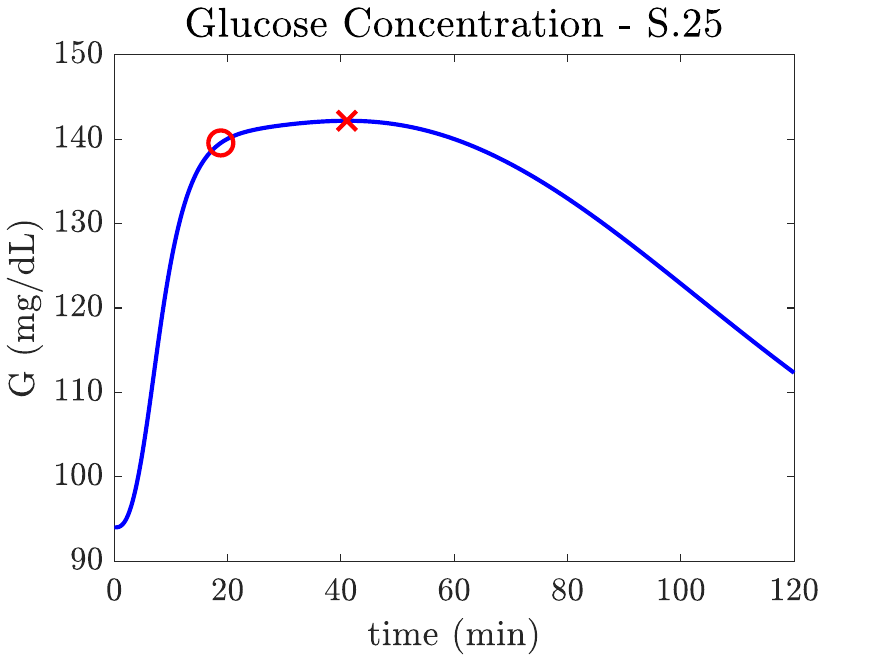}}\\
	\subfloat[]{\includegraphics[width=0.35\linewidth]{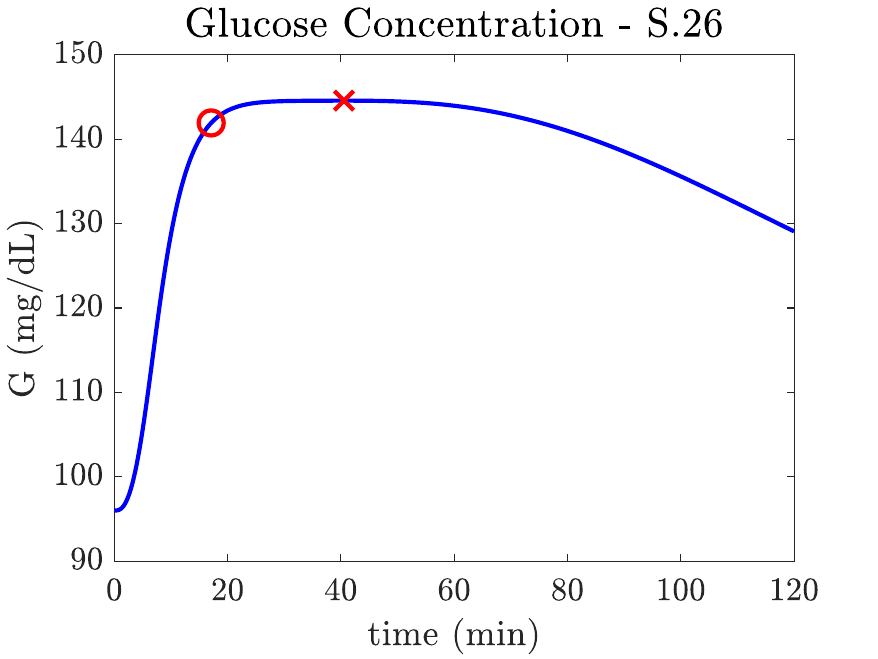}}
	\subfloat[]{\includegraphics[width=0.35\linewidth]{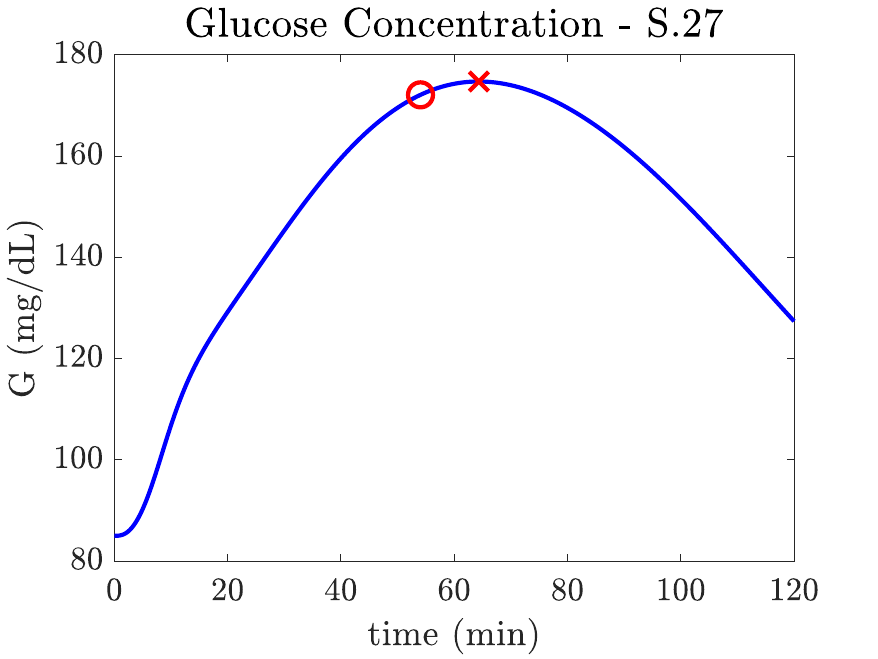}}
	\subfloat[]{\includegraphics[width=0.35\linewidth]{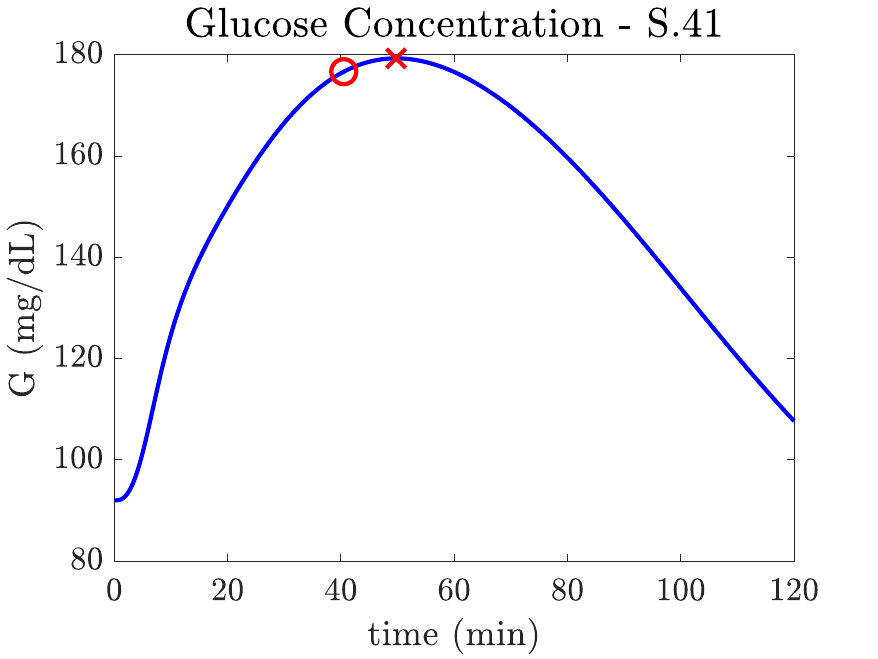}}
	\caption{Glycemic response for the subjects classified as outliers. The red cross indicates the mathematical peak $(\tmax,\Gmax)$, whereas the red circle represents the biological glucose peak $(\tbio,\Gbio)$. For Subjects \textsc{(a)}~12, \textsc{(c)}~25, \textsc{(d)}~26 the biological peak corresponds to the beginning of the plateau.}
	\label{fig:outliers}
\end{figure}

As mentioned before, by looking at Figure~\ref{fig:corr_kabs}, we can identify six outliers (cyan diamonds) since their rate $\Kabs$ suggests a classification into Group 1, but the glycemic peak is realized after 30 min from glucose ingestion. Therefore, their data were excluded from the grouped statistical analyses for further investigation. The glucose response of such subjects is reported in Figure~\ref{fig:outliers}, where a red cross denotes the maximum of $\G$, i.e.  $\Gmax=\max_{t\in[0,120]}{\G(t)}$. 

% CURVE OUTLIERS

We present an automatic procedure to obtain more information from the outliers' responses and understand if they could be classified in one of the identified groups. We distinguish the concept of mathematical peak (that is $\Gmax$) and biological peak $\Gbio$ was distinguished. By biological peak we mean the value of $\G$ which could be considered as the maximum of the glycemic response with respect to the behavior of the glucose-insulin system. For the subjects already classified into the three groups, $\Gbio=\Gmax$.

The biological peak of the outliers can be estimated as $\Gbio=\Gmax-\tolG$, where $\tolG$ is a tolerance and was set to $\tolG=2.6$ mG/dL in the tests. If $\Gbio$ is close to $\Gmax$, then the two quantities basically coincide. On the other hand, if $\Gbio$ is far away from $\Gmax$, it means that the glucose response reaches $\Gbio$ and then shows a very slow increase to $\Gmax$. In this case, we can consider as peak time the instant $\tbio$ associated with $\Gbio$, instead of $\tmax=\mathrm{argmax}_{t\in[0,120]}{\G(t)}$.

We represent $\Gbio$ by a red circle in Figure~\ref{fig:outliers}. It is evident that the glucose response of Subjects 12, 25 and 26 reached a plateau so that the biological peak time $\tbio$ allows classification of those subjects as part of Group~1. On the other hand, $\tbio\approx\tmax$ for Subjects 16, 27 and 41 so that their status as outliers remains unchanged: their biological parameters were similar to those of Group~1, but $\tmax>30$~min so that a clear classification is not possible.

\begin{figure}[t]
	\centering
	\subfloat[]{\includegraphics[width=0.34\linewidth]{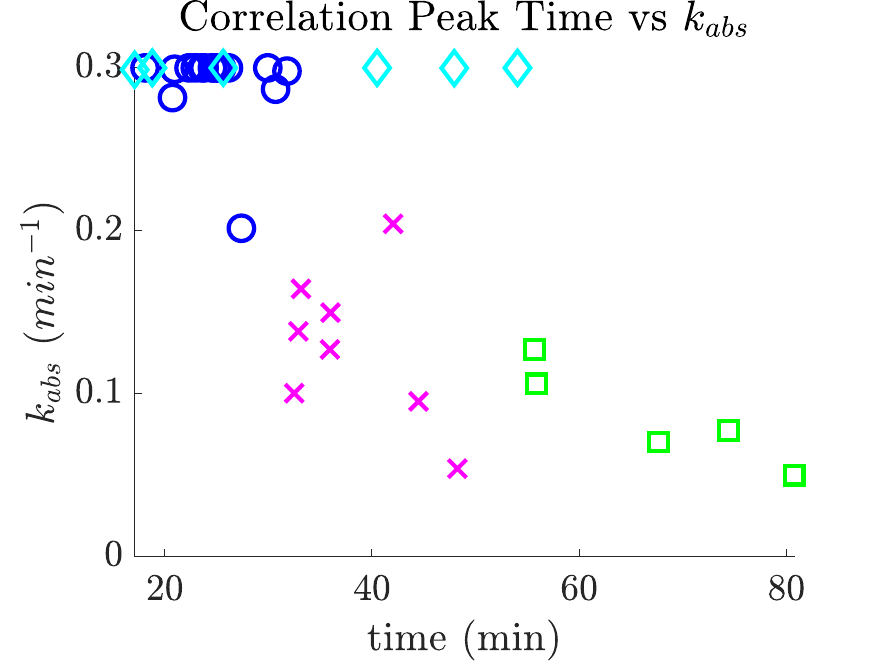}}
	\subfloat[]{\includegraphics[width=0.34\linewidth]{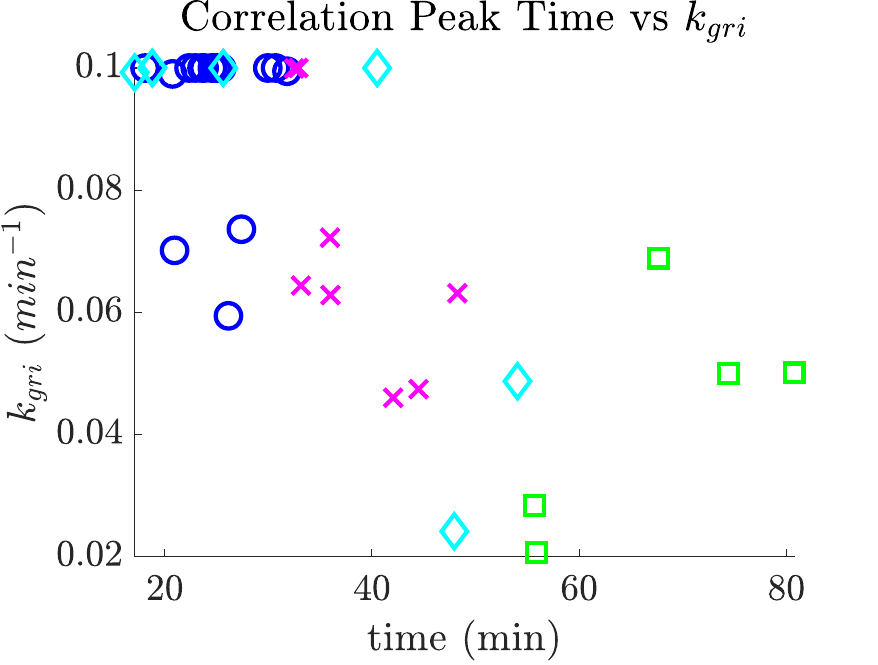}}
	\subfloat[]{\includegraphics[width=0.34\linewidth]{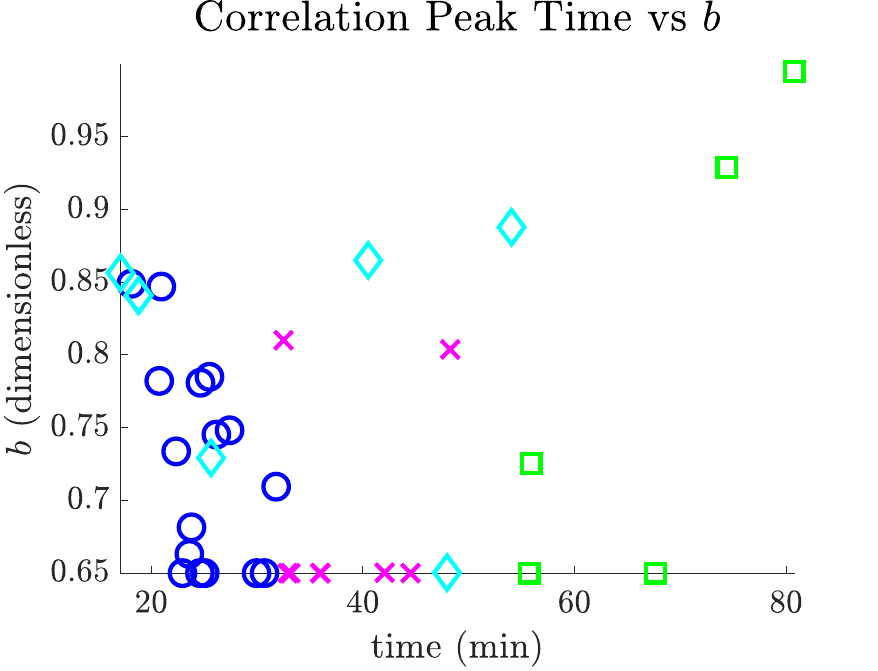}}
	\caption{Correlation plots of \textsc{(a)}~$\Kabs$, \textsc{(b)}~$\Kgri$ and \textsc{(c)}~$\b$ with respect to the time of glycemic peak for Group 1 (blue circles), Group 2 (magenta crosses), Group 3 (green squares). Notice that three of six outliers (cyan diamonds) are now classified as part of Group~1.}
	\label{fig:corr_after_outliers}
\end{figure}

\begin{figure}[t]
	\centering
	\begin{overpic}[width=0.8\linewidth]{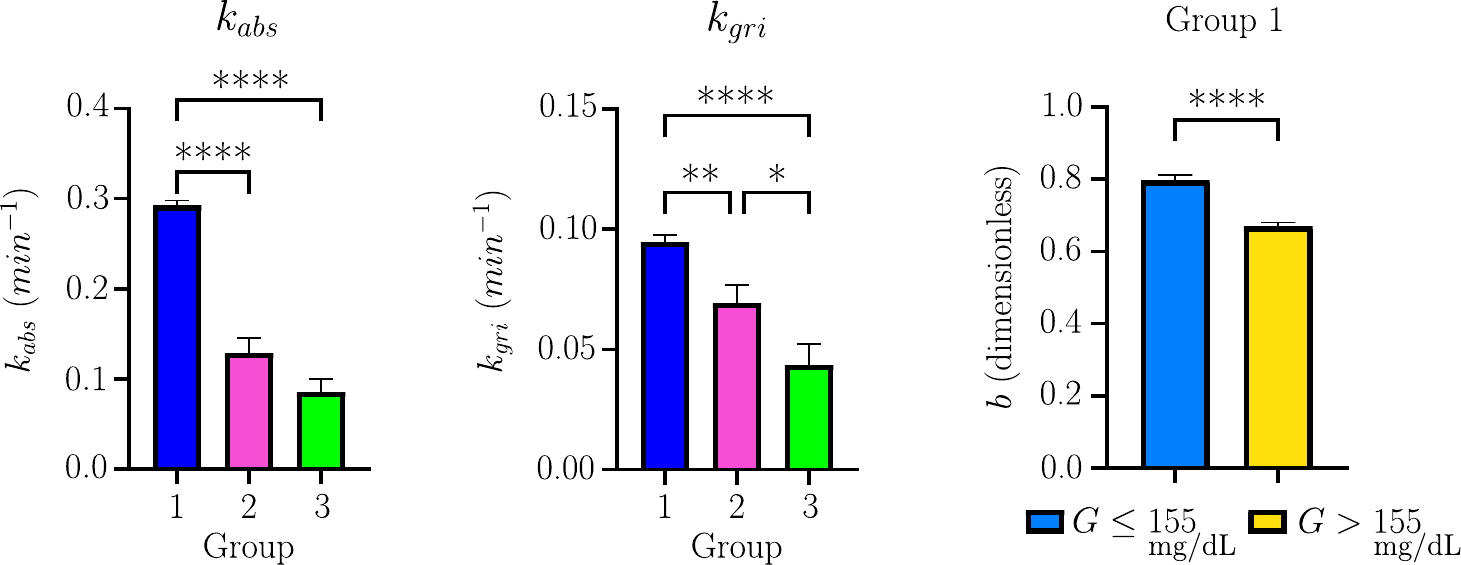}
		\put(52,-12.5){\small\textsc{(a)}}
		\put(166,-12.5){\small\textsc{(b)}}
		\put(280,-12.5){\small\textsc{(c)}}
	\end{overpic}
	
	\
	
	\caption{Histograms of \textsc{(a)}~$\Kabs$, \textsc{(b)}~$\Kgri$ for Group 1 (blue), Group 2 (magenta), Group 3 (green), and \textsc{(c)}~$\b$ of Group 1 for subjects with glucose values $\G\le155$ mg/dL (blue) and $\G>155$ mg/dL (yellow). Statistical significance (One--Way ANOVA followed by Bonferroni \emph{post-hoc} test): $p < 0.05$ (*); $p < 0.01$ (**); $p < 0.001$ (***); $p < 0.0001$ (****).}
	\label{fig:stat_after_outliers}
\end{figure}

After adding Subjects 12, 25 and 26 to Group~1, we recomputed the related statistics; the results are reported in Table~\ref{tab:group1_after_outliers}. Moreover, Figure~\ref{fig:corr_after_outliers} collects the updated correlation plots between the Peak Time and $\Kabs$, $\Kgri$ and $\b$ respectively, while in Figure~\ref{fig:stat_after_outliers} the associated histograms are reported. Specifically, $\Kabs$ is statistically highest in Group 1 compared to Group 2 and Group 3. Similarly, $\Kgri$ was significantly higher in Group 1 than in the other groups. We observe distrinct kinetic behaviors within Group 1, which was further divided based on the value of $\Gbio$, with all comparisons showing statistical significance. 

We emphasize once again that the rate of appearance parameters are an effective tool to analyze glucose responses and identify different metabolic reactions.

% STATISTICHE OUTLIERS

\begin{table}[t]\renewcommand{\arraystretch}{1.2}
	\centering
	\begin{tabular}{ccccc}
		& \# Subjects & Peak Time & $\Kabs$ & $\Kgri$\\
		\hline
		Group 1 & 19 & $24.2026 \pm 0.9315$ & $0.2929 \pm 0.0052$ & $0.0948 \pm 0.0028$\\
		\hline
		\textit{Detail:} & \# Subjects & $\Gbio$ & $\b$ &\\
		\hline
		$\Gbio\le 155$ & 10 & $139.9808\pm2.1676$ & $0.8038\pm0.0149$ & \\
		$\Gbio> 155$ & 9 & $176.8198\pm 3.9325$ & $0.6767 \pm 0.0109$ & \\
		\hline
	\end{tabular}
	\caption{Mean values and standard deviation of Peak Time, $\Kabs$, and $\Kgri$ for Group 1 (after classification of outlier Subjects). In detail, mean values and standard deviation of $\Gbio$ and $\b$ for the sub-groups having $\Gbio\le155$ and $\Gbio>155$. For the subjects we classified at the beginning, $\Gbio=\Gmax$.}
	\label{tab:group1_after_outliers}
\end{table}

%\NOTE{Controllare statistica 4 ****. Aggiungere statistiche nella cpation dei plot istrogrammi}

\section{Conclusion}\label{sec:conclusions}

The standard procedure for computing the Glycemic Index of food is not immediate and easy to perform due to the several requirements listed by the International Standard Organization (ISO): the procedure consists in testing a group of at least ten people, which should adhere to a set of behavioral guidelines, even during the days before the test.

Our final goal is the design of an automatic procedure for computing GI without the need for human testing, and this paper represents a first step in that direction. The idea is to exploit a mathematical model~\cite{dallaman07} to estimate the glucose response, given the composition of the ingested food.

In this work, we focused on the ingestion of pure glucose, whose blood response is the standard reference for computing GI of food. More precisely, we exploited a datasets of \textit{in vivo} measurements \rv{on healthy subjects} to estimate the parameters governing the mathematical model. We observed that the glucose responses provided by the model are well-approximating the laboratory results. In addition, we noticed that the correlation between a subset of parameters (i.e. the rate of appearance parameters describing glucose absorption) and the glucose curves allows us to classify subjects into three groups depending on the time they reach the glycemic peak. 

\rv{Our study presents some limitations concerning the size and type of population. Specifically, our dataset only consists of 35 healthy subjects. In addition, the methodology for the \textit{in vivo} measurements could be improved by performing continuous glucose monitoring, which provides more accurate data compared to the periodic test.}

%This feature may also be employed for further analysis regarding glycaemia and/or disease prevention. 

We plan to extend the results of this preliminary work to larger datasets \rv{with both healthy and diabetic subjects}. Moreover, we will also study the correlation between the parameters governing the mathematical model and the ingestion of composite food, which may contain fat, proteins and fibers beyond glucose/carbohydrates.

%\backmatter

\section*{Declarations}

The authors declare that they have no known competing financial interests or personal relationships that could have appeared to influence the work reported in this paper. 

\

\noindent Data will be made available on request.

\bibliographystyle{abbrv}
\bibliography{biblio}

\end{document}